\documentclass[11 pt]{amsart}
\usepackage{amssymb}
\usepackage{amsmath}
\usepackage{amsfonts}
\usepackage{graphicx}
\usepackage{amsthm}
\usepackage{enumerate}
\usepackage[mathscr]{eucal}
\usepackage{verbatim}
\theoremstyle{plain}
\newtheorem{theorem}{Theorem}[section]
\newtheorem{lemma}[theorem]{Lemma}
\newtheorem{proposition}[theorem]{Proposition}

\numberwithin{equation}{section}
\theoremstyle{plain}

\numberwithin{equation}{section}
\theoremstyle{remark}

\def\intslash{\rlap{\kern  .32em $\mspace {.5mu}\backslash$ }\int}
\def\qsl{{\rlap{\kern  .32em $\mspace {.5mu}\backslash$ }\int_{Q_x}}}
\def\Re{\operatorname{Re\,}}

\def\vthb{\vartheta}

\def\emph#1{{\it #1 }}

\def\inn#1#2{\langle#1,#2\rangle}

\def\noi{\noindent}

\def\lc{\lesssim}
\def\gc{\gtrsim}

\def\eps{\varepsilon}

\def\la{\lambda}

\def\om{\omega}
\def\Om{\Omega}

\def\bbR{{\mathbb {R}}}

\def\cU{{\mathcal {U}}}
\def\cV{{\mathcal {V}}}

 at 10 true pt

\def\be#1{\begin{equation}\label{#1}}
\def\bas{\begin{align*}}
\def\eas{\end{align*}}
\def\bi{\begin{itemize}}
\def\ei{\end{itemize}}

\def\eps{\varepsilon}

\def\emph#1{{\it #1}}
\def\textbf#1{{\bf #1}}
\def\intslash{\rlap{\kern  .32em $\mspace {.5mu}\backslash$ }\int}
\def\qsl{{\rlap{\kern  .32em $\mspace {.5mu}\backslash$ }\int_{Q_x}}}

\begin{document}
\date{Revised version, July 13, 2006}
\title
[Fourier transforms of measures  on curves]
%[Average decay of Fourier transforms]
{Average decay estimates for Fourier transforms of measures supported on curves}

\author[]
%[L. Brandolini, G. Gigante, A. Greenleaf, A. Iosevich, A. Seeger,
%G. Travaglini]
{Luca Brandolini \  \ Giacomo Gigante \ \
Allan Greenleaf\\
Alexander  Iosevich \ \
Andreas Seeger \ \
 Giancarlo Travaglini}
%in $\mathbb{R}^{d}$}
%\author{L. Brandolini\and G. Gigante\and A. Greenleaf
%\and A. Iosevich\and A. Seeger\and G. Travaglini}

\address {L. Brandolini \\ 
Dipartimento di Ingegneria
Gestionale e dell'Informazione
\\ Universit\`a di Bergamo\\  Viale Marconi 5\\
24044 Dalmine BG\\ Italy}

\email{luca.brandolini@unibg.it}

\address {G. Gigante\\
Dipartimento di Ingegneria
Gestionale e dell'Informazione
\\ Universit\`a di Bergamo\\  Viale Marconi 5\\
24044 Dalmine BG\\ Italy}

\email{giacomo.gigante@unibg.it}

\address{A. Greenleaf \\  Department of Mathematics \\
University of Rochester \\
Ro\-ches\-ter, NY 14627, USA}

\email{allan@math.rochester.edu}

\address{A. Iosevich \\
Department of Mathematics\\
University of Missouri-Columbia
\\ Columbia, MO 65211, USA}

\email{iosevich@math.missouri.edu}

\address{A. Seeger   \\
Department of Mathematics\\ University of Wisconsin-Madison\\Madison, WI 53706, USA}
\email{seeger@math.wisc.edu}

\address{G. Travaglini \\
Dipartimento di Matematica e Applicazioni \\
Universit\`a di Milano-Bicocca \\
Via R. Cozzi 53 \\  20125 Milano, Italy}

\email{giancarlo.travaglini@unimib.it}

%\subjclass{42B10, 42B99}
%\keywords{Fourier transforms of measures on curves, $L^q$ average decay,
%Fourier restriction problem}

\thanks{Research supported in part by NSF grants.}

%\thanks{AMS Subject Classification: 42B10, 42B99}
%\thanks{Keywords: Fourier transforms of measures on curves,
%$L^q$ average decay,
%Fourier restriction problem.}

\begin{abstract}
We consider Fourier transforms $\widehat \mu$ of densities supported on
curves in $\bbR^d$. We obtain sharp lower and close to sharp upper bounds
for the decay rates of $\|\widehat \mu(R\cdot)\|_{L^q(S^{d-1})}$, as
$R\to \infty$.
\end{abstract}

\maketitle

\section{Introduction and Statement of Results}

\label{intro}

In this paper we investigate the relation between the geometry of a curve
$\Gamma$ in $\mathbb{R}^{d}$, $d>2,$ and the spherical $L^{q}$ average decay
of the Fourier transform of a smooth density $\mu$ compactly supported on
$\Gamma$.

Let $\Gamma$ be a smooth ($C^{\infty}$) immersed curve in $\mathbb{R}^{d}$
with parametrization $t\rightarrow\gamma(t)$ defined on a compact interval
$I$
%, with $2\leq K\leq d$
and let $\chi\in C^{\infty}$ be supported in the interior of $I$. Let
$\mu\equiv\mu_{\gamma,\chi}$ be defined by
\begin{equation}
\langle\mu,f\rangle=\int f(\gamma(t))\chi(t)dt \label{mudef}%
\end{equation}
and define by $\widehat{\mu}(\xi)=\int\exp(-i\langle\xi,\gamma(t)\rangle
)\chi(t)dt$ its Fourier transform. For a large parameter $R$ we are interested
in the behavior of $\widehat{\mu}(R\omega)$ as a function on the unit sphere,
in particular in the $L^{q}$ norms
\begin{equation}
G_{q}(R)\equiv G_{q}(R;\gamma,\chi):=\Big(  \int|\widehat{\mu}(R\omega
)|^{q}d\omega\Big)  ^{1/q} \label{Gq}%
\end{equation}
where $d\omega$ is the
rotation invariant measure on $S^{d-1}$ induced by Lebesgue
measure in ${\mathbb{R}}^{d}$. The rate of decay depends on the number of
linearly independent derivatives of the parametrization of $\Gamma$. Indeed if
one assumes that for every $t$ the derivatives $\gamma^{\prime}(t)$,
$\gamma^{\prime\prime}(t)$, ..., $\gamma^{(d)}(t)$ are linearly independent
then from the standard van der Corput's lemma (see \cite[page 334]{stein}) one
gets $G_{\infty}(R)=\max_{\omega}|\widehat{\mu}(R\omega)|=O(R^{-1/d})$. If one
merely assumes that at most $d-1$ derivatives are linearly independent then
one cannot in general expect a decay of $G_{\infty}(R)$; one simply considers
curves which lie in a hyperplane. However Marshall \cite{M} showed
that one gets an optimal estimate
for the \textit{$L^{2}$ average decay}, namely
\begin{equation}
G_{2}(R)=O(R^{-1/2}) \label{generalL2bd}%
\end{equation}
as $R\rightarrow\infty$,
 for every
compactly supported $C^{1}$ curve $\gamma$.

We are interested in estimates for the $L^{q}$ average decay, for
$2<q<\infty$. If $\gamma$ is a straight line such extensions fail,
and additional
conditions are necessary. Our first result addresses the case of nonvanishing
curvature.
%is an extension to values of $q>2$ under
%the  assumption, with a possible logarithmic blowup.

\begin{theorem}
\label{upperI} Suppose that for all $t\in I$ the vectors $\gamma^{\prime}(t)$
and $\gamma^{\prime\prime}(t) $ are linearly independent. Then for $R\ge2$

(i)
\begin{equation}
G_{q}(R)\lesssim%
\begin{cases}
R^{-1/2}\left(  \log R\right)  ^{1/2-1/q}\quad & \text{ if }2\leq q\leq4\\
R^{-2/q}\left(  \log R\right)  ^{1/q}\quad & \text{ if }4\leq q\leq\infty.
\end{cases}
\label{logest}%
\end{equation}

(ii) Suppose that there is $N\in{\mathbb{N}}$ so that for every $\omega\in
S^{d-1}$ the function $s\mapsto\langle\omega,\gamma^{\prime\prime}(s)\rangle$
changes sign at most $N$ times on $I.$ Then
\begin{equation}
\label{qle4sharp}G_{q}(R) \lesssim R^{-1/2} \quad\text{ if } 2\le q< 4
\end{equation}
and
\begin{equation}
\label{qge4sharp}G_{q}(R) \lesssim R^{-2/q} \quad\text{ if } 4< q\le\infty.
\end{equation}

\end{theorem}

Here and elsewhere the notation $a\lesssim b$ means $a\le C b$ for a suitable
nonnegative constant $C$.
%We remark that for the proof of Theorem \ref{upperI}
%the $C^\infty$ smoothness assumption on $\gamma$
%can be replaced by a $C^2$ smoothness assumption.

The $L^{4}$ estimate $G_{4}(R)=O(R^{-1/2}[\log R]^{1/4})$ is sharp even for
nondegenerate curves, \textit{cf.} Theorem \ref{lowerthm} below. The estimate
\eqref{qle4sharp} is sharp and it is open whether for $q\neq4$ there exists an
example for which the logarithmic term in \eqref{logest} is necessary.

The estimate \eqref{qge4sharp} is sharp in the case where the curve lies in a
two dimensional subspace. Under stronger nondegeneracy assumptions this
estimate can be improved. In particular one is interested in the case of
nondegenerate curves in ${\mathbb{R}}^{d}$, meaning that for all $t$ the
vectors $\gamma^{(j)}(t)$, $j=1,\dots, d$, are linearly independent. In the
case $d=2$ we have of course the optimal bound $G_{q}(R)=O(R^{-1/2})$ for all
$q\le\infty$, by the well known stationary phase bound (for results for
general curves in $\bbR^2$ and hypersurfaces in higher dimensions see
\cite{BHI} and  references contained  therein). The situation is more
complicated for nondegenerate curves in higher dimensions, and Marshall
\cite{M} proved (essentially) optimal results for nondegenerate curves in
${\mathbb{R}}^{d}$ if $d=3$ and $d=4$.

We show that one gets close to optimal results for nondegenerate curves in all
dimensions. Our method is different from the explicit computations in
Marshall's paper and relies on a variable coefficient analogue of the Fourier
restriction theorem due to Fefferman and Stein in two dimensions, see
\cite{F}, and due to Drury \cite{D}
for curves in higher
dimensions. The variable coefficient analogues are
due to Carleson and Sj\"olin
\cite{CS} (see also H\"ormander \cite{H}) in two dimensions and to Bak and Lee
\cite{BL} in higher dimensions.

To formulate our result let, for $1\leq q\leq\infty,$%
\begin{equation}
\sigma_{K}(q)\equiv\sigma_{K}^{d}(q)=\left\{
\begin{array}
[c]{ll}%
\min_{\left\{  k=2,\,\ldots,d\right\}  }\frac{1}{k}+\frac{k^{2}-k-2}{2kq}, &
\mathrm{for\,\,}K=d,\\
\min_{\left\{  k=2,\,\ldots,K\right\}  }\left\{  \frac{1}{k}+\frac{k^{2}%
-k-2}{2kq},\,\frac{K}{q}\right\}  , & \mathrm{for\,\,}2\leq K<d.
\end{array}
\right.  \label{sigmadef}%
\end{equation}

\begin{theorem}
\label{sigmathm} Suppose that for all $t\in I$ the vectors $\gamma^{\prime
}(t)$, ..., $\gamma^{(K)}(t) $ are linearly independent. Then for $R\ge2$
\begin{equation}
\label{sigmabd}G_{q}(R)\le C_{\sigma}R^{-\sigma}, \quad\sigma<\sigma_{K}(q).
\end{equation}

\end{theorem}

 For integers $k\ge 1$ set
\begin{equation}
\label{qkdef}q_{k}:= \frac{k^{2}+k+2}2
\end{equation}
so that $q_1=2$,  $q_2=4$, $q_3=7$, $q_4=11$. 
Observe that the set of points $(q^{-1},\,\mathbb{\sigma}_{d}^{d}(q))$, $q\ge 2$,  is the
broken line joining the points
%\[
%(1,\tfrac{1}{2}),\,(\tfrac{1}{4},\tfrac{1}{2}),\,(\tfrac{1}{7},\tfrac{3}%
%{7}),\ldots,(\tfrac{2}{k^{2}+k+2},\tfrac{2k}{k^{2}+k+2}),\ldots,
%(\tfrac
%{2}{d^{2}-d+2},\tfrac{2d-2}{d^{2}-d+2}),\,(0,\tfrac{1}{d}),
%\]
%%%%%%%%%%%%%%%
%\[
%(1,1/2),\, (q_2^{-1},{2}{q_2^{-1}}),
%\ldots,
%(q_k^{-1},{k}{q_k^{-1}}),\ldots,
%(q_{d-1}^{-1},{(d-1)}{q_{d-1}^{-1}}),\,(0,d^{-1}).
%\]
%%%%%%%%%%%%%
\[
(\frac{1}{q_1},\frac{1}{q_1}),
%(\frac {1}{2},\frac{1}{2}),
(\frac{1}{q_2},\frac{2}{q_2}),
\ldots,
(\frac{1}{q_k},\frac{k}{q_k}),\ldots,
(\frac{1}{q_{d-1}},\frac{d-1}{q_{d-1}}),(0,\frac{1}{d}),
\]
while for $K<d,$ the set of points $(q^{-1},\,\mathbb{\sigma}_{K}^{d}(q))$ is
the concave broken line joining the points
\[
%(\frac {1}{2},\frac{1}{2}),\,
(\frac{1}{q_1},\frac{1}{q_1}), (\frac{1}{q_2},\frac{2}{q_2}),
\ldots,(\frac{1}{q_k},\frac{k}{q_k}),\ldots,
(\frac{1}{q_K},\frac{K}{q_{K}}),(0,0).
\]
%%%%%%%%%%%%%%%%%
%\[
%(1,1/2),\, (q_2^{-1},{2}{q_2^{-1}}),
%\ldots,
%(q_k^{-1},{k}{q_k^{-1}}),\ldots,
%(q_K^{-1},{K}{q_{K}^{-1}}),\,(0,0).
%\]
%%%%%%%%%%%%%%%%%
%\[
%(1,\tfrac{1}{2}),\,(\tfrac{1}{4},\tfrac{1}{2}),\,(\tfrac{1}{7},\tfrac{3}%
%{7}),\ldots,(\tfrac{2}{k^{2}+k+2},\tfrac{2k}{k^{2}+k+2}),\ldots,(\tfrac
%{2}{K^{2}+K+2},\,\tfrac{2K}{K^{2}+K+2}),(0,\,0).
%\]
Furthermore observe that $\sigma_{K}^{d}(q)>2/q$ if $3\leq K\leq d$, and 
$q>4$.
The picture shows the  graph $\{1/q,\sigma_{K}^{K}(q)\}$ as a function of $1/q$, for $K=10$.

%TCIMACRO{\FRAME{dhFU}{4.8784in}{3.3235in}{0pt}{\Qcb{The graph of $\sigma
%_{10}^{10}\left(  1/q\right)  .$}}{}{sigmadieci.eps}%
%{\special{ language "Scientific Word";  type "GRAPHIC";
%maintain-aspect-ratio TRUE;  display "USEDEF";  valid_file "F";
%width 4.8784in;  height 3.3235in;  depth 0pt;  original-width 6.5683in;
%original-height 4.4633in;  cropleft "0";  croptop "1";  cropright "1";
%cropbottom "0";  filename 'sigmadieci.eps';file-properties "XNPEU";}}}%
%BeginExpansion
\begin{center}
\includegraphics[
height=3.3235in,
width=4.8784in
]%
{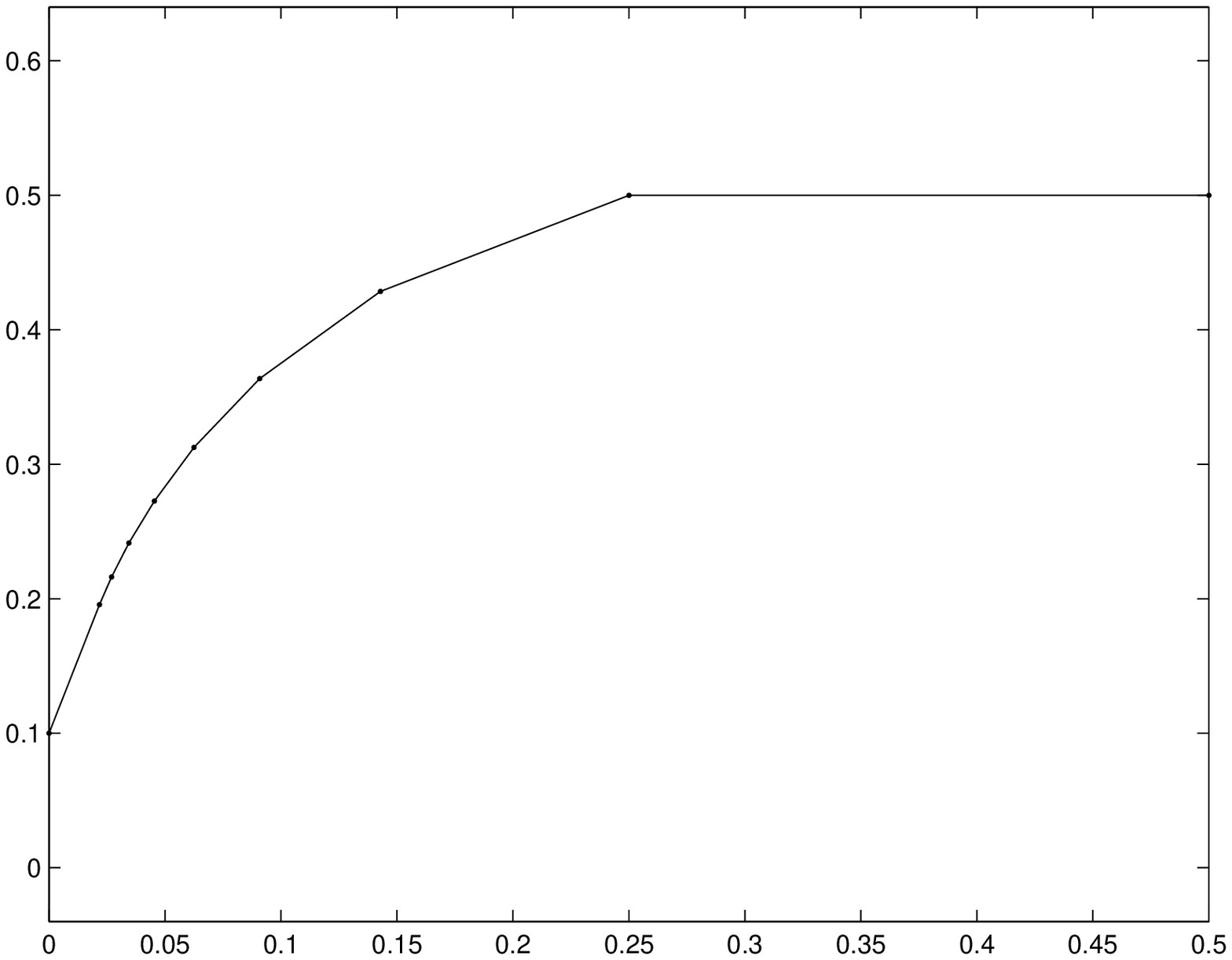}%
\\
The graph  $\{1/q,\sigma_{10}^{10}(q)\}.$%
%The graph of $\sigma_{10}^{10}\left(  1/q\right)  .$%
\end{center}
%EndExpansion
%%%%%%%
%
%\begin{align*}
%& \left(  1,\,\frac{1}{2}\right)  ,\,\left(  \frac{1}{4},\,\frac{1}{2}\right)
%,\,\left(  \frac{1}{7},\,\frac{3}{7}\right)  ,\ldots,\,\left(  \frac{2}%
%{k^{2}+k+2},\,\frac{2k}{k^{2}+k+2}\right)  ,\ldots\\
%& \ldots,\,\left(  \frac{2}{K^{2}+K+2},\,\frac{2K}{K^{2}+K+2}\right)
%,\,\left(  0,\,0\right)  .
%\end{align*}
%
%
%
%
%%%%%%%
%TCIMACRO{\FRAME{dtbpFU}{5.0593in}{3.4477in}{0pt}
%{\Qcb{The graph of $\sigma
%_{6}^{10}\left(  1/q\right)  .$}}{}{sigmacinque.eps}%
%{\special{ language "Scientific Word";  type "GRAPHIC";
%maintain-aspect-ratio TRUE;  display "USEDEF";  valid_file "F";
%width 5.0593in;  height 3.4477in;  depth 0pt;  original-width 6.7366in;
%original-height 4.5777in;  cropleft "0";  croptop "1";  cropright "1";
%cropbottom "0";  filename 'sigmacinque.eps';file-properties "XNPEU";}}}%
%BeginExpansion
%%%%%%%
%
%%%SECOND PICTURE, slightly too large
%\begin{center}
%\includegraphics[
%height=3.4477in,
%width=5.0593in
%]%
%{sigmacinque.eps}%
%\\
%The graph  $\{1/q,\sigma_{6}^{10}(q)\}.$%
%%The graph of $\sigma_{6}^{10}\left(  1/q\right)  .$%
%\end{center}
%%EndExpansion
%%%%%%%
 We emphasize that the graph
of $\sigma^{d}_{K}$ is slightly different for $d>K$, as then the left 
line segment connects 
$(q_K^{-1},Kq_K^{-1})$ to $(0,0)$.

Theorem \ref{sigmathm} is sharp up only to endpoints, at least for
nondegenerate curves (for which $\gamma^{\prime}(t)$,..., $\gamma^{(d)}(t)$
are linearly independent), and also for some other cases where $K<d$,
$\gamma^{\prime}(t)$,..., $\gamma^{(K)}(t)$ are independent and $\gamma$ lies
in a $K$ dimensional affine subspace. We note that for the case $K=d=4$
Marshall \cite{M} obtained the sharp bound $G_{q}(R)\lesssim R^{-\sigma
_{4}^{4}(q)}$ when $4<q<7$ and $q>7$; moreover $G_{q}(R)\lesssim
R^{-\sigma_{4}^{4}(q)}\log^{1/q}(R)$ when $q=4$ or $q=7$ (the logarithmic term
for the $L^{4}$ bound seems to have been overlooked in \cite{M}).

\medskip

We now state lower bounds for the average decay.
The cutoff function $\chi$ is as in \eqref{mudef} (and $G_q(R)$ depends on
$\chi$).

\begin{theorem}
\label{lowerthm} Suppose that $2\le K\le d$ and, for some $t_{0}\in I$, the
vectors $\gamma^{\prime}(t_{0})$, ..., $\gamma^{(K)}(t_{0}) $ are linearly
independent. Then for suitable $\chi\in C^\infty_0$ there are
 $c>0$,
$R_0\ge 1$ so that the
following lower
bounds  hold for $R>R_0$.

(i) If $2\le K\le d-1$ then
%we have the lower bounds%
\begin{equation}
\label{withoutlog}G_{q}(R)\ge c R^{-\sigma_{K}(q)}, \quad2<q< q_{K};
\end{equation}
moreover
\begin{equation}
\label{withlog}G_{q}(R)\ge c R^{-\sigma_{K}(q)} \log^{1/q}(R), \quad
q\in\{q_{k}: k=2,\dots, K-1\}.
\end{equation}

(ii) If $K=d$ then
\begin{align}
\label{lowerdwolog}
G_{q}(R)  &  \ge c R^{-\sigma_{d}(q)}, \quad2<q\le\infty,\\
\label{lowerdwlog}
G_{q}(R)  &  \ge c R^{-\sigma_{d}(q)} \log^{1/q}(R), \quad q\in\{q_{k}:
k=2,\dots, d-1\}.
\end{align}

(iii) If $2\le K\le d-1$ and, in addition, $\gamma^{(K+1)}\equiv0$ then
\begin{align}
\label{lowerKwolog}
G_{q}(R)  &  \ge c R^{-\sigma_{K}(q)}, \quad2<q\le\infty,\\
\label{lowerKwlog}
G_{q}(R)  &  \ge c R^{-\sigma_{K}(q)} \log^{1/q}(R), \quad q\in\{q_{k}:
k=2,\dots, K\}.
\end{align}

\end{theorem}
%In the case of the Veronese curves,
%\[\gamma_{K}(t)=\sum_{j=1}^{K} t^{j} e_{j},\quad2\le K\le d,\]
%these results say that
%\[G_{q}(R)\gtrsim R^{-\sigma_{K}(q)},\]
%and\[G_{q_{k}}(R)\gtrsim R^{-\sigma_{K}(q_{k})}\log^{1/q_{k}}(R),
%\]
%for $2\le k\le K$ if $K<d$ and for $2\le k<d$ if $K=d$.

%\begin{proposition} \label{sigmalowerthree}
%Suppose  that for some  $t_0\in I$
%the vectors  $\gamma'(t_0)$, $\gamma''(t_0)$ and $\gamma'''(t_0)$ are
%linearly independent.
%Then   $\chi$ in \eqref{mudef} can be chosen so that
%for sufficiently large
%$R$
%\begin{equation}\label{thirdsigmalowerbd}
%G_4(R)\ge C R^{-1/2} [\log R]^{1/4}
%\end{equation}
%\end{proposition}

\medskip

\noi\textit{Remark:}
A careful examination of the proof  yields
some  uniformity in the lower bound.
Assume that $\gamma^{(j)}(t_0)=e_j$,
(the $j$th unit vector), $j=1,\dots, K$, and
$\|\gamma\|_{C^{K+3}(I)}\le C_1$. Then there is $h=h(C_1)>0$
so that for every  smooth $\chi$ supported in $(-h,h)$ with
$\Re \chi(t)>c_1>0$ in $(-h/2,h/2)$ there exists an $R_0$
depending only on
$c_1$, $C_1$, $\|\chi'\|_\infty$ and $\|\chi''\|_\infty$
 so that the above lower bounds hold for $R\ge R_0$. We shall
not pursue  this point in detail.

\medskip

\noi{\bf Addendum:}
 After the first version of this paper had been submitted
we learned about the
work of
Arkhipov, Chubarikov and Karatsuba \cite{akt}, \cite{akt2} who proved
sharp estimates for the
$L^q(\Bbb R^d)$ norms  of the Fourier transform of smooth
densities on certain polynomial curves.
We are grateful to Jong-Guk Bak  who pointed out these references to us.
 The work of these authors  shows that
for, say $\gamma(t)=\sum_{k=1}^d t^ke_k$, $t\in[0,1]$ the Fourier transform
$\widehat {d\sigma}$ belongs to $L^q(\bbR^d)$ if and only if $q>q_d=
(d^2+d+2)/2$.
This  result seems to have been overlooked until recently; it
rules out an $L^{q_d'}$
 endpoint bound  for the Fourier  restriction problem associated to
curves, {\it cf.}  a discussion
in \cite{mock} and a remark in \cite{BL}).
More can be said  in two dimensions where  the
endpoint restricted weak type (4/3) inequality for the 
Fourier restriction 
operator is known to fail by a Kakeya set argument, see \cite{bcss}.
We note that the
 lower bound in chapter 2 of
\cite{akt2} is  closely related to
\eqref{lowerdwolog}
and the method in \cite{akt2} actually can be used to yield
\eqref{lowerdwolog} for the curve $(t,...,t^d)$ in the range
$q\ge q_{d-1}$; vice versa  one notices $\sigma_d(q_d)=d/q_d$
and   integrates the lower bound for
$R^{d-1}G_{q}^{q}(R)$ in $R$ to obtain lower bounds for
$\|\widehat {d\sigma}\|_{L^{q}(\bbR^d)}$.

A variant of an argument in \cite{akt2} can be shown to close the $\eps$ gap
between  upper and lower bounds in some cases. We formulate one such result.

\begin{theorem}\label{fttype}
Suppose that $\gamma$ is smooth and is  either of finite type, or  polynomial.

Assume that $\gamma'(t),\dots, \gamma^{(K)}(t)$ are  linearly independent,
for every $t\in I$. Then the following holds:
%%Let  $q_k=(k^2+k+2)/2$. Then

(i) If $K=d$ then
\begin{equation}
\label{withoutlogu}
G_{q}(R)\le C_q R^{-\sigma_{d}(q)}, \quad
q\ge 2, \quad
q\notin\{q_{k}: k=2,\dots, d-1\},
\end{equation}
and
\begin{equation}
\label{withlogu}
G_{q}(R)\lc R^{-\sigma_{d}(q)} \log^{1/q}(R), \quad
q\in\{q_{k}: k=2,\dots, d-1\}.
\end{equation}

(ii) If $ 2\le K\le d-1$ then
\begin{equation}
\label{withoutloguK}
G_{q}(R)\le C_q R^{-\sigma_{K}(q)}, \quad
q\ge 2, \quad
q\notin\{q_{k}: k=2,\dots, K\}.
\end{equation}
and
\begin{equation}
\label{withloguK}
G_{q}(R)\lc R^{-\sigma_{K}(q)} \log^{1/q}(R), \quad
q\in\{q_{k}: k=2,\dots, K\}.
\end{equation}
\end{theorem}

It is understood that the implicit constants in 
\eqref{withoutlogu} and 
\eqref{withoutloguK} depend on $q$ as $q\to q_k$.
Note that in the finite type case \eqref{withoutloguK} and 
 \eqref{withloguK} can be improved for $q>q_{K}$ since we have  
some nontrivial decay for $G_\infty(R)$. However, for the 
sharpness in the most 
degenerate case compare Theorem \ref{lowerthm}, part (iii).

\medskip

\noi{\it Remark:} The result of Theorem
\ref {fttype}, for polynomial curves,  could be used to obtain the upper bounds of
Theorem \ref{sigmathm}, which involves a loss of $R^\eps$, by a polynomial
approximation argument.
Note however, that such an argument  requires upper bounds
for derivatives of $\gamma$
 up to order $C+\eps^{-1}$, as $\eps\to 0$. An examination of the proof of
Theorem \ref{sigmathm} shows that one can get away with 
upper  bounds for the 
derivatives up to order $N$ where $N$ depends on the dimension but
not on $\eps$.

\medskip

\noindent\textit{Structure of the paper.} In \S \ref{upI} we prove the
estimates \eqref{qle4sharp} and \eqref{qge4sharp} which involve the assumption
of $\langle\omega,\gamma^{\prime\prime}(s)\rangle$ not changing sign. Here we
also discuss an application to some mixed norm inequalities for rotated
measures. In \S \ref{upII} we prove Theorem \ref{sigmathm} and
\eqref{logest}.
In \S\ref{upIII} we give the proof of Theorem \ref{fttype}.
%In \S\ref{upIII} we prove
%an additional result for curves in $\bbR^d$ which is related to
%Marshall's $\bbR^4$ bound.
%%%removed  this section.
In \S \ref{revisited} we revisit some known asymptotic expansion with precise
quantifications which are convenient for the proof of the lower bounds. The
proof of the lower bounds of Theorem \ref{lowerthm} is given in
\S \ref{frombelow}.

\section{Upper bounds, I.}
\label{upI}

We shall now prove part (ii) of Theorem \ref{upperI}
(\textit{i.e.} \eqref{qle4sharp}, \eqref{qge4sharp})
 under  the less
restrictive smoothness
condition $\gamma\in C^{2}(I)$;
we recall the assumptions that $\gamma^{\prime}(t)$
and $\gamma^{\prime\prime}(t)$ are linearly
independent and that
we also require
that the functions
$s\mapsto \langle\omega,\gamma^{\prime\prime}(s)\rangle$ have at most a
bounded number of
sign changes
on $I$. Note that this hypothesis is certainly satisfied
 if $\gamma$ is a polynomial, or a trigonometric polynomial, or smooth and
of finite type.

%We shall now prove Theorem \ref{upperI},
%under  a less restrictive smoothness
%condition:
%\begin{theorem}\label{upper2I} Suppose that $\gamma \in C^2(I)$ and that
%the vectors $\gamma^{\prime}(t)$
%and $\gamma^{\prime\prime}(t) $ are linearly independent,
%for all $t\in I$.
% Then
%inequalities \eqref{logest} hold. If in addition
%the functions $s\mapsto\langle\omega,\gamma^{\prime\prime}(s)\rangle$,
%$\omega\in S^{d-1}$ have at most a bounded number of sign changes
%then \eqref{qle4sharp} and \eqref{qge4sharp} hold.\end{theorem}
% Note that the second  hypothesis on sign changes is certainly satisfied
% if $\gamma$ is a polynomial, or a trigonometric polynomial,
%or smooth and
%of finite type.

We need a result on oscillatory integrals which is a consequence of the
standard van der Corput Lemma; it is also related to a more sophisticated
statement on oscillatory integrals with polynomial phases in \cite{Ob}.

Let $\eta$ be a $C^{\infty}$ function with support in $(-1,1)$ so that
$\eta(s)=1$ in $(-1/2,1/2)$; we also assume that $\eta^{\prime}$ has only
finitely many sign changes. Let $\eta_{1}(s)=\eta(s)-\eta(2s)$ (so that
$1/4\leq|s|\leq1$ on $\operatorname*{supp}\eta_{1}$) and let
\[
\eta_{l}(s)=\eta_{1}(2^{l-1}s)
\]
so that $2^{-l-1}\leq\left\vert s\right\vert \leq2^{-l+1}$ on the support of
$\eta_{l}$.

\begin{lemma}
\label{oilemma} Let $I$ be a compact interval and let $\chi\in C^{1}(I)$. Let
$\phi\in C^{2}(I)$ and suppose that $\phi^{\prime\prime}$ changes signs at
most $N$ times in $I$.

Then, for $1\leq2^{l}\leq\lambda$,
\[
\left\vert \int_{I}\eta_{l}(\phi^{\prime}(s))e^{i\lambda\phi(s)}%
\chi(s)ds\right\vert \leq CN2^{l}\lambda^{-1}.
\]

\end{lemma}

\begin{proof}
We may decompose $I$ into subintervals $J_{i}$, $1\leq i\leq K$, $K\leq2N+2$,
so that both $\phi^{\prime}$ and $\phi^{\prime\prime}$ do not change sign in
each $J_{i}$. Each interval $J_{i}$ can be further decomposed into a bounded
number of intervals $J_{i,k}$ so that $\eta^{\prime}\left(  \phi^{\prime
}\right)  $ is of constant sign in $J_{i,k}$. It suffices to estimate the
integral ${\mathcal{I}}_{i,k}$ over $J_{i,k}$. By the standard van der Corput
Lemma, the bound ${\mathcal{I}}_{i,k}=O(2^{l}/\lambda)$ follows if we can show
that
\[
\int_{J_{i,k}}\left\vert \partial_{s}\left(  \eta_{l}(\phi^{\prime}%
(s))\chi(s)\right)  \right\vert ds\leq C
\]
which immediately follows from
\begin{equation}
\int_{J_{i,k}}\Big\vert 2^{l}\phi^{\prime\prime}(s)\eta_{1}^{\prime}
(2^{l}\phi^{\prime}(s))\Big\vert ds\leq C.
\end{equation}
But by our assumption on the signs of $\phi^{\prime}$, $\phi^{\prime\prime}$,
and $\eta^{\prime}$ the left hand side is equal to
\[
\Big\vert \int_{J_{i,k}}2^{l}\phi^{\prime\prime}(s)\eta_{1}^{\prime}
(2^{l}\phi^{\prime}(s))ds\Big\vert =\Big\vert \int_{J_{i,k}}\partial
_{s}\Big(  \eta_{l}(\phi^{\prime}(s))\Big)  ds\Big\vert \leq C.
\]

\end{proof}

\begin{proof}
[Proof of \eqref{qle4sharp} and \eqref{qge4sharp}] We may assume that $\Gamma$
is parametrized by arc\-length and that the support of $\chi$ is small (of
diameter $\ll1$). Determine the integer $M(R)$ by $2^{M}\leq R<2^{M+1}.$ With
$\eta_{l}$ as above
%Moreover set $\eta^M(s)=\eta(2^Ms)$
%which is supported on $[-2^{-M}, 2^{-M}]$.
define for $l<M$
\[
g_{R,l}(\omega)=\int e^{iR\langle\omega,\gamma(s)\rangle}\chi(s)\eta_{l}
(\langle\omega,\gamma^{\prime}(s)\rangle)ds
\]
and for $l=M$ define $g_{R,M}$ similarly by replacing the cutoff $\eta
_{l}(\langle\omega,\gamma^{\prime}(s)\rangle)$ with $\eta(2^{M} \langle
\omega,\gamma^{\prime}(s)\rangle)$. We can decompose
\[
\int e^{iR\langle\omega,\gamma( t)\rangle}\chi( t) dt=\sum_{l\leq M}%
g_{R,l}(\omega)
\]
and observe that $g_{R,l}=0$ if $l\leq-C$.

It follows from Lemma \ref{oilemma} that
\begin{equation}
\sup_{\omega\in S^{d-1}}|g_{R,l}(\omega)|\lesssim2^{l}/R \label{infty}%
\end{equation}
We also claim that
\begin{equation}
\Big ( \int|g_{R,l}(\omega)|^{2}d\omega\Big)^{1/2}\lesssim%
\begin{cases}
2^{l}R^{-1} & \text{ if }2^{l}\leq R^{1/2},\\
2^{-l}(1+\log(2^{2l}R^{-1}))^{1/2} & \text{ if }2^{l}\geq R^{1/2}.
\end{cases}
\label{two}%
\end{equation}
Given (\ref{infty}) and (\ref{two}) we deduce that
\begin{align}
\left\Vert g_{R,l}\right\Vert _{L^{q}(S^{d-1})}  &  \leq\left\Vert
g_{R,l}\right\Vert _{L^{2}(S^{d-1})}^{2/q}\left\Vert g_{R,l}\right\Vert
_{L^{\infty}(S^{d-1})}^{1-2/q}\nonumber\\
&  \lesssim%
\begin{cases}
2^{l}R^{-1} & \text{ if }2^{l}\leq R^{1/2},\\
2^{l(1-4/q)}(1+\log(2^{2l}R^{-1}))^{1/q}R^{-1+2/q} & \text{ if }2^{l}\geq
R^{1/2}.
\end{cases}
\label{Lqinequalityforfixedl}%
\end{align}
If $q\neq4$ the asserted bound $O(R^{-2/q})$ bound follows by summing in $l$.

We now turn to the proof of (\ref{two}). Note that (\ref{two}) follows
immediately from (\ref{infty}) if $2^{l}\leq R^{1/2}$. Now let $2^{l}\geq
R^{1/2}$. For the $L^{2}$ estimate in this range we shall just use the
nonvanishing curvature assumption on $\Gamma$. We need to estimate the $L^{2}$
norm of $g_{R,l}$ over a small coordinate patch ${\mathcal{V}}$ on the sphere
where we use a regular parametrization $y\rightarrow\omega(y)$, $y\in
\lbrack-1,1]^{d-1}$; \textit{i.e.}
\begin{equation}
\left\vert \int u(y)g_{R,l}(\omega(y))\overline{g_{R,l}(\omega(y))}
dy\right\vert \lesssim2^{-2l}(1+\log(2^{2l}R^{-1})) \label{lm},
\end{equation}
where $u\in C_{0}^{\infty}$, so that $\omega(y)\in{\mathcal{V}}$ if
$y\in\operatorname*{supp}(u)$. The left hand side of (\ref{lm}) can be written
as
\begin{multline*}
{\mathcal{I}}_{l}:=\iint_{s_{1},s_{2}}\int_{y}u(y)e^{iR \langle\omega
(y),\gamma(s_{1})-\gamma(s_{2})\rangle}\chi(s_{1}) \overline{\chi(s_{2})}\\
\times\eta(2^{l}\langle\omega(y),\gamma^{\prime}(s_{1})\rangle)\eta(2^{l}
\langle\omega(y),\gamma^{\prime}(s_{2})\rangle)dyds_{1}ds_{2}%
\end{multline*}
and we note that on the support of the amplitude we get that $\gamma^{\prime
}(s_{i})$ is almost perpendicular to $\omega$, \textit{i.e.} we may assume by
the assumption of small supports that there is a direction $w=(w_{1}%
,\dots,w_{d-1})$ so that
\[
\Big|\sum_{\nu=1}^{d-1}w_{\nu}\partial_{y_{\nu}}\langle\omega(y),\gamma
^{\prime}(s)\rangle\Big| \geq1/2
\]
if $s\in\operatorname*{supp}(\chi)$ and $y\in\operatorname*{supp}(u)$. By a
rotation in parameter space we may assume that
\begin{equation}
\big| \partial_{y_{1}} \langle\omega(y),\gamma^{\prime}(s)\rangle
\big| \geq1/2\quad\text{ if }s\in\operatorname*{supp}(\chi),\quad
y\in\operatorname*{supp}(u). \label{lowerbdforyderivatives}%
\end{equation}

Now let for fixed unit vectors $v_{1}$, $v_{2}$ and $\delta>0$
\[
{\mathcal{U}}_{\delta}(v_{1},v_{2})=\{\omega\in S^{d-1}:| \langle\omega
,v_{1}\rangle|\leq\delta,|\langle\omega, v_{2}\rangle|\leq\delta\}
\]
and observe that the spherical measure of this region is at most $O(\delta)$;
moreover this bound can be improved if $|v_{1}-v_{2}|$ is $\geq\delta$. Namely
if $\alpha(v_{1},v_{2})$ is the acute angle between $v_{1}$ and $v_{2}$ then
\begin{equation}
{\text{\textrm{meas}}}({\mathcal{U}}_{\delta}(v_{1},v_{2}))\lesssim
\min\left\{  \delta,\frac{\delta^{2}}{\sin\alpha(v_{1},v_{2})}\right\}
\label{Udeltaestimate}%
\end{equation}

The condition (\ref{lowerbdforyderivatives}) implies that
\[
\left\vert \langle\partial_{y_{1}}\omega(y), \gamma(s_{1})-\gamma
(s_{2})\rangle\right\vert \geq c|s_{1}-s_{2}|
\]
and given the regularity of the amplitude we can gain by a multiple
integration by parts in $y_{1}$ provided that $|s_{1}-s_{2}|\geq2^{l}/R$;
indeed we gain a factor of $O(R^{-1}|s_{1}-s_{2}|^{-1}2^{l})$ with each
integration by parts. We obtain, for any $N$,
\begin{align}
{\mathcal{I}}_{l}  &  \lesssim\int_{s_{1}\in\operatorname*{supp}(\chi)}\left[
\int_{|s_{2}-s_{1}|\leq c}\operatorname*{meas}\left(  {\mathcal{U}}_{2^{-l}%
}(\gamma^{\prime}(s_{1}),\gamma^{\prime}(s_{2}))\right)  \right.
\label{part}\\
&  \qquad\qquad\qquad\left.  \times\min\{1,(R|s_{1}-s_{2}|2^{-l}%
)^{-N}\}^{^{^{^{^{{}}}}}}ds_{2}\right]  ds_{1}.\nonumber
\end{align}

By the assumption that $|\gamma^{\prime\prime}(s)|$ is bounded below and
$\gamma^{\prime}$ and $\gamma^{\prime\prime}$ are orthogonal we get as a
consequence of \eqref{Udeltaestimate}
\[
\operatorname*{meas}\left(  {\mathcal{U}}_{2^{-l}}(\gamma^{\prime}%
(s_{1}),\gamma^{\prime}(s_{2}))\right)  \leq\min\{2^{-l},2^{-2l}|s_{1}%
-s_{2}|^{-1}\}.
\]
Now we use this bound and integrate out the $s_{2}$ integral in (\ref{part})
and see that the main contribution comes from the region where $2^{-l}
\leq|s_{1}-s_{2}|\leq2^{l}/R$ which yields the factor $\log(R2^{-2l})$ in
(\ref{lm}).
\end{proof}

\medskip

\noindent\textbf{An application.} We consider a $C^{2}$ curve $\gamma:\left[
-a,a\right]  \rightarrow\mathbb{R}^{d}$ with nonvanishing curvature and assume
that, as in \eqref{qle4sharp}, the function $s\mapsto
\langle\omega,\gamma^{\prime\prime}(s)\rangle$ has a bounded number of sign changes.

Let $\mu$ be the measure induced by the Lebesgue measure on $\Gamma$,
multiplied by a smooth cutoff function. For every $\sigma\in SO\left(
d\right)  $ define $\mu_{\sigma}$ by $\mu_{\sigma}\left(  E\right)
=\mu\left(  \sigma E\right)  $ and for every test function $f $ in
$\mathbb{R}^{d}$%
\[
Tf\left(  x,\sigma\right)  =f\ast\mu_{\sigma}\left(  x\right)  .
\]
We are interested in the $L^{p}({\mathbb{R}}^{d})\to L^{s}(SO(d),
L^{q}({\mathbb{R}}^{d}))$ mapping properties, in particular for $q=p^{\prime
}=p/(p-1)$. This question had been investigated in \cite{RT} for curves in the
plane, with essentially sharp results in this case, see also \cite{BGT}. The
standard example, namely testing $T$ on characteristic functions of balls of
small radius yields the necessary condition $1+(d-1)/q\ge d/p$. Setting
$q=p^{\prime}$ we see that the $L^{p}({\mathbb{R}}^{d})\to L^{s}(SO(d),
L^{p^{\prime}}({\mathbb{R}}^{d}))$ fails for $p<(2d-1)/d$ (independent of $s$).

The approach in \cite{RT} together with the inequality \eqref{qle4sharp}
yields
\begin{equation}
\| Tf\|_{L^{s}(SO(d),L^{p^{\prime}}({\mathbb{R}}^{d}))}\leq
C_{p}\Vert f\Vert_{p},\quad p=\tfrac{2d-1}{d},\quad s<\tfrac{4d-2}{d}.
\label{mixednorm}%
\end{equation}

\begin{proof}
[Proof of \eqref{mixednorm}]We imbed $T$ in an analytic family of operators.
After rotation and reparametrization (modfying the cutoff function) we may
assume that $\gamma(t)=\sum_{j=1}^{d-1}\varphi_{j}(t)e_{j}+te_{d}$, with
$\varphi_{j}(0)=0$. Let $z\in\mathbb{C}$ such that $\operatorname{Re}z>0$ and
define a distribution $i_{z}$ by $\left\langle i_{z},\chi\right\rangle
=(\Gamma(z))^{-1}\int_{0}^{+\infty}\chi\left(  t\right)  t^{z-1}dt$. Then
define $\mu_{\sigma}^{z}$ by $\widehat{\mu_{\sigma}^{z}}(\xi)=\widehat
{\mu_{\sigma}}(\xi)\prod_{j=1}^{d}\widehat{i_{z}}(\langle\sigma\xi
,e_{j}\rangle)$ and $T^{z}$ by $T^{z}f(x,\sigma)=\mu_{\sigma}^{z}\ast f$.
Following \cite{RT} one observes that $\mu_{\sigma}^{1+i\lambda}$ is a bounded
function, namely we have
\[
|\langle\mu^{1+i\lambda},g\rangle|\lesssim\int_{\mathbb{R}}\int_{{\mathbb{R}%
}^{d-1}}|g(x_{d}e_{d}+\sum_{j=1}^{d-1}(y_{j}+\phi_{j}(x_{d})e_{j}%
)|dy_{1}\cdots dy_{d-1}dx_{d}\lesssim\Vert g\Vert_{1}%
\]
so that
\begin{equation}
T^{1+i\lambda}:L^{1}(\mathbb{R}^{d})\rightarrow L^{\infty}(SO(d)\times
\mathbb{R}^{d}). \label{Bound2}%
\end{equation}
We also have
\begin{equation}
T^{-\tfrac{1}{2d-2}+i\lambda}:L^{2}(\mathbb{R}^{d})\rightarrow L^{q}%
(SO(d),L^{2}(\mathbb{R}^{d})),\quad2\leq q<4. \label{Bound1}%
\end{equation}
The implicit constants in both inequalities are at most exponential in
$\lambda$. Thus we obtain the assertion \eqref{mixednorm} by analytic
interpolation of operators.

To see \eqref{Bound1} we observe that $\widehat{i_{z}}(\tau)=O(|\tau
|^{-\operatorname{Re}(z)})$ and apply Plancherel's theorem and then
Minkowski's integral inequality to bound for $\alpha>0$
\begin{align*}
\big\|&T^{-\alpha+i\lambda}f\Vert_{L^{q}(L^{2})}^{2}  \\& =
\Big(\int_{SO(d)}\Big(
\int_{{\mathbb{R}}^{d}}|\widehat{f}(\xi)\widehat{\mu}_{\sigma}|^{2}\prod
_{j=1}^{d-1}|\widehat{i_{-\alpha+i\lambda}}(\langle\sigma\xi,e_{j}
\rangle)|^{2}d\xi\Big)  ^{q/2}d\sigma\Big)^{2/q}\\
&  \lesssim\int_{\mathbb{R}^{d}}\Big\vert \widehat{f}(\xi)\Big\vert
^{2}\Big(  \int_{SO(d)}\Big\vert \widehat{\mu_{\sigma}}(\xi)\Big\vert
^{q}\prod_{j=1}^{d-1}|\langle\sigma\xi,e_{j}\rangle|^{\alpha q}d\sigma
\Big)
^{2/q}d\xi,
\end{align*}
and by \eqref{qle4sharp} and the assumption $q<4$ the last expression is
dominated by a constant times
\[
\int_{\mathbb{R}^{d}}\left\vert \widehat{f}(\xi)\right\vert ^{2}|\xi
|^{2\alpha(d-1)}\Big(\int_{SO(d)}\big|\widehat{\mu_{\sigma}}(\xi
)\big|^{q}d\sigma\Big)^{2/q}d\xi\lesssim\int_{\mathbb{R}^{d}}\big|\widehat
{f}(\xi)\big|^{2}|\xi|^{2\alpha(d-1)-1}d\xi.
\]
For $\alpha=(2d-2)^{-1}$, this yields the bound \eqref{Bound1}.
\end{proof}

\noindent\textit{Remark.} We do not know whether the index $s=\frac{4d-2}%
{d-1}$ in \eqref{mixednorm} is sharp. The following example only shows that we
need $s\leqslant10$ for $d=3$. Let $\gamma\left(  t\right)  =\left(
t,t^{2},0\right)  $ and let $\chi_{B_{\delta}}$ be a box centered at the
origin with sides parallel to the axes and having sidelengths $1,1$ and
$\delta$. A computation shows that $| T\chi_{B_{\delta}}( x,\sigma) |
\geqslant c$ for $\sigma$ in a set of measure $\varepsilon^{2}$ and $x$ in a
set of measure $\varepsilon$, for some small $\varepsilon>0$. It follows that $p^{-1}\le2s^{-1}+q^{-1}$. For
$p=5/3$ and thus $p^{\prime}=5/2$ this yields $s\le10$.

\section{Upper bounds, II}

\label{upII} We are now concerned with the proof of Theorem \ref{sigmathm} and
the proof of part (i) of Theorem \ref{upperI}. For the latter
%In order to prove part (i) of Theorem
%\ref{upperI}
we use a version of the Carleson-Sj\"olin theorem (\cite{CS}, \cite{H}), and
for Theorem \ref{sigmathm} we use a recent generalization due to Bak and Lee
\cite{BL}. These we now recall.

Consider, for large positive $R$,
\[
T_{R}f\left(  x\right)  =\int_{\mathbb{R}}e^{iR\phi\left(
x,\,t\right)  }a\left(  x,\,t\right)  f\left(  t\right)  dt
\]
with real valued phase function $\phi\in C^{\infty}({\mathbb{R}}^{n}%
\times{\mathbb{R}})$, and compactly supported amplitude $a\in{C}_{0}^{\infty}(
\mathbb{R}^{n}\times{\mathbb{R}}).$ Assume the non-vanishing torsion
condition
\begin{equation}
\det\left(  \partial_{t}\left(  \nabla_{x}\phi\right)  ,\,\partial_{t}%
^{2}\left(  \nabla_{x}\phi\right)  ,\,\ldots,\partial_{t}^{n}\left(
\nabla_{x}\phi\right)  \right)  \neq0 \label{carlesonsjolin}%
\end{equation}
on the support of $a.$ Then if $p^{-1}+n(n+1)(2q)^{-1}=1$ and $q>(n^{2}%
+n+2)/2$ there is a constant $C_q$ independent of $f$ and of $R\geq2$ such that
\begin{equation}
\|T_{R}f\|_{L^{q}( \mathbb{R}^{n})}\leq C_qR^{-n/q}\|f\|_{L^{p}%
(\mathbb{R}) }. \label{baklee}%
\end{equation}
When $n=2$ it is well known that a slight modification of H\"{o}rmander's
proof (\cite{H}) of the Carleson-Sj\"{o}lin theorem gives the 
endpoint result
\begin{equation}
\left\vert \left\vert T_{R}f\right\vert \right\vert _{L^{4}\left(
\mathbb{R}^{n}\right)  }\lc R^{-1/2}\log^{1/4}R\left\vert \left\vert
f\right\vert \right\vert _{L^{4}\left(  \mathbb{R}\right)  };
\label{hormander}%
\end{equation}
see also \cite{MSe}, where a somewhat harder vector-valued analogue is proved.

In order to establish estimates \eqref{logest} we need to show that under the
assumption of linear independence of $\gamma^{\prime}(t)$ and $\gamma
^{\prime\prime}(t)$ (for each $t\in I$) that
\begin{equation}
G_{4}(R)\lesssim R^{-1/2}[\log R]^{1/4}. \label{elfour}%
\end{equation}
To establish \eqref{sigmabd} under the assumption that the first $K$
derivatives are linearly independent for every $t\in I$ we need to show that
for any $2\le k\le K$
\begin{equation}
G_{q}(R)\lesssim R^{-k/q}, \quad q>q_{k}, \label{elq}%
\end{equation}
where $q_{k}$ is as in \eqref{qkdef}. All other estimates in \eqref{logest},
\eqref{sigmabd} follow by the usual convexity property of the $L^{p}$ norm
(\textit{i.e.} $\|F\|_{p}\le\| F\|_{p_{0}}^{1-\vartheta} \| F\|_{p_{1}%
}^{\vartheta}$ for $p^{-1}=(1-\vartheta)p_{0}^{-1}+\vartheta p_{1}^{-1}$).

\begin{proof}
[Proof of \eqref {elq} and \eqref {elfour}]Let
\begin{equation}
\label{FRdefinition}
F_{R}(\omega)= \int_{\mathbb{R}}e^{i\langle\omega,\gamma(t)\rangle}a ( \omega)
\chi( t) dt.
\end{equation}
By compactness, we can suppose that $\chi$ is supported in $(-\varepsilon
,\varepsilon) $, with $\varepsilon$ as small as we need. Divide the sphere
$S^{d-1}$ into two subsets $A$ and $B$; here in $A,$ the unit normal to the
sphere is essentially orthogonal to the span of the vectors $\gamma^{\prime
}(0) ,\ldots,\gamma^{(k)}(0) ,$ and in $B,$ the unit normal to the sphere is
close to the span of $\gamma^{\prime}(0) ,\ldots,\gamma^{( k) }(0)$.
%Assume
%$q>\left(  k^{2}+k+2\right)  /2$.

Now consider a coordinate patch $V$ of diameter $\varepsilon$ on $A$ and
parametrize it by $y\mapsto\omega( y) $ with $y\in\mathbb{R}^{d-1},$ $y$ near
$Y_{0}.$ From the defining property of $A,$ it follows that the vectors
$\nabla_{y} \langle\omega(\cdot),\gamma^{(j)}(t)\rangle$, $j=1,\dots, k$ are
linearly independent when evaluated at $y$ near $Y_{0}$, provided that
$|t|<\varepsilon$. Therefore we can choose the parameterization $y=(x^{\prime
},y^{\prime\prime})=(x_{1},\ldots,x_{k},y_{k+1},\ldots,y_{d-1})$ in such a way
that also the vectors $\nabla_{x^{\prime}} \langle\omega(\cdot),
\gamma^{(j)}(t)\rangle$, $j=1,\dots, k$ are linearly independent.
%the projections of $\gamma^{\prime}\left(  0\right)  ,\ldots
%,\gamma^{\left(  k\right)  }\left(  0\right)  $
%on the hyperplane tangent to
%$S^{d-1}$ in $\omega\left(  Y_{0}\right)  $ generate a parallelepiped of
%$k$-dimensional volume greater than some $\delta>0.$ If we choose the $y$
%coordinates in such a way that
%$\partial_{y_{1}}\omega\left(  Y_{0}\right)
%,\ldots,\,\partial_{y_{k}}\omega\left(  Y_{0}\right)  $
%form an orthonormal
%system whose span contains the parallelepiped,
%then we may express its volume
%as
%\begin{equation}
%\left\vert \det\left(
%\begin{array}
%[c]{ccc}%
%\partial_{y_{1}}\omega\left(  Y_{0}\right)  \cdot\gamma^{\prime}\left(
%$0\right)  & \cdots & \partial_{y_{1}}\omega\left(  Y_{0}\right)  \cdot
%\gamma^{\left(  k\right)  }\left(  0\right) \\
%\vdots & \ddots & \vdots\\
%\partial_{y_{k}}\omega\left(  Y_{0}\right)  \cdot\gamma^{\prime}\left(
%0\right)  & \cdots & \partial_{y_{k}}\omega\left(  Y_{0}\right)  \cdot
%\gamma^{\left(  k\right)  }\left(  0\right)
%\end{array}
%\right)  \right\vert >\delta>0.\label{orthogonality}%
%\end{equation}
If we consider $y^{\prime\prime}$ as a parameter and we define
\[
\phi^{y^{\prime\prime}}( x^{\prime},\,t) = \langle\omega(x^{\prime}%
,y^{\prime\prime}),\gamma( t)\rangle,
\]
then the phase functions $\phi^{y^{\prime\prime}}$ satisfy condition
(\ref{carlesonsjolin}) uniformly in $y^{\prime\prime}$. We also have upper
bounds for the higher derivatives in $\omega$ and $\gamma$ which are
uniform
in $y^{\prime\prime}$ as well (here $y^{\prime\prime}$ is taken from a
relevant compact set). Thus one can apply the Bak-Lee result \eqref{baklee} in $k$
dimensions to obtain, for fixed $y^{\prime\prime}$,
\begin{equation}
\label{Aestimate}\Big(\int\Big| \int_{\mathbb{R}}e^{iR\phi^{y^{\prime\prime}}(
y^{\prime},\,t) }a ( \omega( y^{\prime},y^{\prime\prime}) ) \chi( t) dt
\Big|^{q} dy^{\prime}\Big)^{1/q} \lesssim R^{-k/q}, \quad q>q_{k}.
\end{equation}
An integration in $y^{\prime\prime}$ yields $\|F_{R}\|_{L^{q}(V)} \lesssim
R^{-k/q}$ for $q>( k^{2}+k+2) /2.$ Similarly if $k=2$ and $q=4$ we can apply
\eqref{hormander} in two variables to obtain $\|F_{R}\|_{L^{4}(V)} \lesssim
R^{-1/2} \log^{1/4}R$. This settles the main estimate for the $L^{q}(A)$ norm.
As for contribution of the $L^{q}(B)$ norm we recall that the unit normal to
the sphere is close to the span of $\gamma^{\prime}(0) ,\ldots, \gamma^{( k)
}( 0) ,$ and thus $\sum_{j=1}^{k}|\langle\omega,\gamma^{(j)}(t)\rangle|>0$.
Therefore we can apply van der Corput's lemma and obtain the $L^{\infty}$
estimate%
\begin{equation}
\label{Linfinity}\|F_{R}\|_{L^{\infty}(B)}\lesssim R^{-1/k}%
\end{equation}
For $k=2,$ this completes the proof of the theorem. For $2<k\leq K,$ we argue
by induction. We assume that 
the asserted estimate holds for $k-1,$ ($k\ge 3$);
that is
\begin{equation}
\label{induction-h}\|F_{R}\|_{q} \lc R^{-(k-1)/{q}} \quad\text{ for }
q>q_{k-1} =\tfrac{k^{2}-k+2}{2}
\end{equation}
where the implicit constants depend on $q$. 
If $\vartheta_{k}=1-q_{k-1}/q_{k}$ then we use the relation $q_{k}-q_{k-1}=k$
to verify that $(1-\vartheta_{k})(k-1)/q_{k-1}+\vartheta_{k}/k=k/q_{k}$. Thus
by a convexity argument we see that a combination of \eqref{Linfinity} and
\eqref{induction-h} yields that
\[
\|F_{R}\|_{L^{q}(B)} \lesssim R^{-k/q}, \quad\text{for $q>\tfrac{k^{2}%
+k+2}2=q_{k}$}.
\]
Together with the corresponding bound for $\|F_{R}\|_{L^{q}(A)}$ proved
above, this concludes the proof.
\end{proof}

\section{Upper bounds, III}\label{upIII}
We give the proof of Theorem \ref{fttype}
 under the {\it finite type assumption.}
By compactness, there is an integer $L\ge d$ and a constant $c>0$ so that
%We say that $\gamma $ is of  finite type if
%there is $L\ge 1$ and $c>0$ so that
for every $s\in I$ and every  $\theta\in S^{d-1}$ we have $ \sum_{n=1}^{L} |\inn{\gamma^{(n)}(s)}{\theta}|\ge c.$
%me that there is $N$ so that for every $\theta\in S^{d-1}$ the function
%\inn{\gamma^{(K)}(t)}{\theta}$ changes sign not more than $N$ times.

We shall argue by induction on $k$.
 By Theorem
\ref{sigmathm}
 the conclusion holds for $k=2.$
Assume $k>2$, and that the desired inequalities are  already proved for
$2\le q\le q_{k-1}$.

Let $F_R$ be as in
\eqref{FRdefinition} and assume that the cutoff function $\chi$ is supported in $(-\eps,\eps)$.
As in the proof of Theorem
\ref{sigmathm}
we split the sphere into subsets $A$ and $B$
where in  $A,$ the unit normal to the
sphere is almost  perpendicular  to the span of the vectors $\gamma^{\prime}(t) ,\ldots,\gamma^{(k)}(t) ,$  for all $|t|<\eps$ 
and in $B,$ the projections of the  unit normals to the sphere
 to  the span of $\gamma^{\prime}(t) ,\ldots,\gamma^{( k) }(t)$
have length $\ge c>0$.

We shall estimate the $L^{q}(A\cap \Omega)$ norm of $F_R$ on a small  patch $\Om$  on the
sphere, and by further localization
we may assume by the finite type assumption 
that there is an  $n\le L$ so that
\begin{equation} \label{n-localiz}
|\inn   {\gamma^{(n)}(s)}{\theta}|\ge c>0, |s|\le \eps, \theta\in \Om.
\end{equation}

We distinguish between the case  $n\ge k$ and $n<k$.
First we assume $n\ge k$ (the main case).
Then  there is the pointwise bound
\begin{equation} \label{HR}
F_R(\theta)\lc \min\big\{1, H_R(\theta))^{-1}\}
\end{equation}
where
$$
H_R(\theta)=\min_{s\in I}
\max_{1\le j\le n}R^{1/j}| \inn{\gamma^{(j)}(s)}{\theta}|^{1/j}.
$$
This is immediate from van der Corput's   lemma; indeed the finite type 
assumption
allows  the decomposition of  the interval $[-\eps,\eps]$ into a
bounded number of subintervals
so that on each subinterval  all derivatives
of $s\mapsto
\inn{\gamma^{(j)}(s)}{\theta}$, $1\le j\le n-1$  are monotone and one-signed.
We now have to estimate the $L^q(A\cap\Omega)$ norm of the right hand side of
\eqref{HR}.

For an $l>0$ consider the set
\begin{equation}\label{Omlr}
\Omega_l(R)=\{\theta\in \Om: H_R(\theta)\in [2^{l}, 2^{l+1})\}.
\end{equation}
By \eqref{n-localiz} we have
$|H_R(\theta)|\gc R^{1/n}.$
Thus
only the values with
\begin{equation} \label{lrestriction}
2^l\gc R^{1/n}
\end{equation}
are relevant
(and likewise the set
of
$\theta\in \Om$ for which $ H_R(\theta)\lc 1$ is empty if $R$ is large).

By the definition of $H_R$ we can find a point  $s_*=s_*(\theta)$ and an integer
 $j_*$, $1\le j^*\le n$,
so that
$H_R(\theta)= |R\inn{\gamma^{(j_*)}(s_*)}{\theta}|^{1/j_*}$
and
$|R\inn{\gamma^{(j)}(s_*)}{\theta}|^{1/j} \le H_R(\theta)$
for all $\theta\in \Om$ and all $j\le n$. This implies
\begin{equation}\label{sstarbd}
|\inn{\gamma^{(j)}(s_*)}{\theta}| \lc 2^{(l+1)j } R^{-1}, \text{ if }\theta
\in \Om_l(R), \quad j\le n.
\end{equation}

We shall now apply a nice idea of \cite{akt}: We divide our interval
 $(-\eps,\eps)$ into $O(2^l)$ intervals $I_{\nu,l}$ of length $\approx 2^{-l}$,
with right endpoints $t_\nu$, so that $t_{\nu}-t_{\nu-1}\approx 2^{-l}$.
The point $s_*$ lies in one of these intervals, say in  $I_{\nu_*}$.
We  estimate
$|\inn{\gamma^{(j)}(t_{\nu_*})}{\theta}|^{1/j}$ in terms of
$H_R(\theta)$.
By a Taylor expansion
we get
\begin{equation} \label{taylorexpansion}
\inn{\gamma^{(j)}(t_{\nu_*})}{\theta}=
\sum_{r=0}^{n-j-1}
\inn{\gamma^{(j+r)}(t_{\nu_*})}{\theta}\frac{(t_{\nu_*}-s_*)^r}{r!}
+
\inn{\gamma^{(n)}(\widetilde t)}{\theta}\frac{(t_{\nu_*}-s_*)^{n-j}}{(n-j)!}
\end{equation}
where $\widetilde t $ is between $s_*$ and $t_{\nu_*}$.
By  \eqref{sstarbd} the  terms in the sum are all  $O(2^{lj}R^{-1})$.
The remainder  term is  $O(2^{-l(n-j)})$
which is also
$O(2^{lj}R^{-1})$, by the condition
\eqref{lrestriction}.
%$2^l\ge c R^{1/n}$.
Now define
%Thus if for some large constant $C$ we define
$$\Om_{\nu,l}=\{\theta\in \Omega:
|\inn{\gamma^{(j)}(t_\nu)}{\theta}|\le C 2^{lj}R^{-1},
j=1,\dots, n\}
$$
and if $C$ is sufficiently large
then
the set $\Omega_l(R)$ is contained in the union of the sets
$\Omega_{\nu,l}$; the constant $C$ can be chosen independently of $l$ and $R$.
%%and we have $O(2^{l})$ such sets).

%As $n\ge k$ the set
%$$\Om_{\nu,l}^*=\{\theta\in \Omega:
%|\inn{\gamma^{(j)}(t_\nu)}{\theta}|\le C 2^{lj}R^{-1},
%j=1,\dots, k\}
%%$$
%contains $\Om_{\nu,l}$.
In view of the linear independence of the vectors
$\gamma^{(j)}(t_\nu)$, $j=1,\dots,k$ and the condition
$\theta\in A$, the measure of the set
$\Om_{\nu,l}$
is $O(\prod_{s=1}^k (2^{sl}R^{-1}))=O(2^{l k(k+1)/2} R^{-k})$, for every
$1\le \nu\lc 2^l$, and thus the measure of the  set
$\Omega_l(R)$ is $O(2^{l(k^2+k+2)/2} R^{-k})$.
On $\Omega_l(R)$ we have
$|F_R(\theta)|\le H_R(\theta)^{-1}\lc 2^{-l}$. Therefore
\begin{equation}\label{kderiv}
\int_{\Om\cap A}|F_R(\theta)|^q d\theta \lc \sum_{cR^{1/n}\le 2^l\le c R}
2^{-lq}
2^{l(k^2+k+2)/2} R^{-k},
\end{equation}
which yields the endpoint bound
$$
\Big(\int_{\Om\cap A}|F_R(\theta)|^{q_k} d\theta\Big)^{1/q_k} \lc R^{-k/q_k}
(\log R)^{1/q_k}.
$$
Of course we also get (by using the same argument with just $k-1$ derivatives)
\begin{equation}\label{k-1deriv}
\int_\Om|F_R(\theta)|^q d\theta \lc \sum_{cR^{1/n}\le 2^l\le c R}
2^{-lq}
2^{l(k^2-k+2)/2} R^{1-k}
\end{equation}
which yields the sharp $L^{q_{k-1}}(A\cap \Omega)$ bound.
Now we consider $q$ satisfying 
$q_{k-1}<q<q_k$,  $k<d$  or  $q_{d-1}<q<\infty$ and  $K=d$.
 We distinguish the cases
(i)  $2^l< R^{1/k}$ and (ii)  $2^l\ge R^{1/k}$.
In the first case we use
\eqref{kderiv} while in the second case we use
\eqref{k-1deriv}.
Then in the case $k<d$
\begin{multline*}
\Big(\int_\Om|F_R(\theta)|^{q} d\theta\Big)^{1/q}
\\
\lc 
R^{-k}
\Big(\sum_{2^l< R^{1/k}}
2^{l(-2q+k^2+k+2)/2} +\sum_{2^l\ge  R^{1/k}}
2^{l(-2q+k^2-k+2)/2} R
\Big)^{1/q}
\end{multline*}
which is bounded by $C R^{-1/k-(k^2-k-2)/2kq}$ if $q_{k-1}<q<q_k$.
If $K=k=d$ then only values with $2^l\ge R^{1/d}$ are relevant and 
only the second sum in the last displayed line occurs. Thus if $K=d$ we
obtain the  estimate
$C R^{-1/d-(d^2-d-2)/2dq}$ for $q>q_{d-1}$.

Now if $n<k$ one gets even better bounds;  we use the induction hypothesis.
First note that for $n=1,2$ integration by parts,
 or van der Corput's lemma,
yields a better bound; therefore assume $n\ge 3$.
We have the bounds
%\begin{align*}
%&\|F_R\|_{L^\infty(\Om)}\lc R^{-1/n}\\
%&\|F_R\|_{L^{q_n}(\Om)}  \lc R^{-n/q_n}(\log R)^{1/q_n};\end{align*}
$\|F_R\|_{L^\infty(\Om)}\lc R^{-1/n}$
and
$\|F_R\|_{L^{q_n}(\Om)}  \lc R^{-n/q_n} (\log R)^{1/q_n}$;
the first one by van der Corput's lemma and the second one by the
induction hypothesis. By convexity this yields
$\|F_R\|_{L^{q_k}}\lc R^{-\alpha(k,n)}\log R^{1/q_k}$
where $\alpha(k,n)= n/q_k+(1-q_n/q_k)/n$ and one checks
that $\alpha(k,n)= k/q_k+ (k-n)(k+1-n)/(2nq_k)$ if   $n<k$,
so that one gets a better estimate.
 The case
 $q_{k-1}<q<q_k$, $n<k$  is handled in the same way.
This yields the desired bounds for the $L^q(A)$ norm of $F_R$.

For the  $L^q(B)$ bound we may use van der Corput's estimate
with $\le k$ derivatives to get an $L^\infty$ bound $O(R^{-1/k})$; we interpolate this with the appropriate $L^p$ bound for $q_{k-2}<p\le q_{k-1}$ which
holds by the induction hypothesis;
 the argument is similar to that in the proof of Theorem \ref{sigmathm}.
This finishes the argument  under the finite type assumption.

\medskip

\noi{\it  Modification for polynomial curves:}
If the coordinate functions  $\gamma_j$ are polynomials of degree $\le L$
we need to take $n=L$ in the definition of $H_R(\theta)$.
We use, for the case $l>0$, the analogue of the
  Taylor expansion
\eqref{taylorexpansion} up to order $L$ with zero remainder term
(again $n=L$).
As above we obtain for $l>0$
the bound
\[
\int_{\Om_l(R)}|F_R(\theta)|^q d\theta \lc  2^{-lq}  R^{-k}
\min\{
2^{l(k^2+k+2)/2}, 2^{l(k^2-k+2)/2} R\}.
\]
Summing in $l>0$ works as before. However we also have  a
contribution from the set
$\Omega_0(R)=\{\theta\in \Om: H_R(\theta)\le 1\}.$
By the polynomial assumption  a  Taylor expansion
(now about the point $s_*$, without remainder) is used to show that
 $\Omega_0(R)$ is contained in the subset of $A$
where $|\inn{\gamma^{(j)}(s_*)}{\theta}|\le CR^{-1}$, $j=1,\dots, k$.
This set has measure $O(R^{-k})$. Thus the desired
bound for $l=0$ follows as well.\qed

\section{Asymptotics for oscillatory integrals revisited}

\label{revisited}

We examine the behavior of some known asymptotics for oscillatory integrals
under small perturbations. This will be used in the subsequent section to
prove the lower bounds of Theorem \ref{lowerthm}.

For $k=2,3,\dots$, there is the following formula for $\lambda>0$:
\begin{equation}
\label{Jla}\int_{-\infty}^{\infty}e^{i\lambda s^{k}} ds= \alpha_{k}
\lambda^{-1/k},
\end{equation}
where
\begin{equation}
\label{ak}\alpha_{k}=%
\begin{cases}
\tfrac2k\Gamma(\tfrac1k)\sin(\tfrac{(k-1)\pi}{2k}), \, & k \text{ odd},\\
\tfrac2k\Gamma(\tfrac1k) \exp(i\tfrac{\pi}{2k}), \, & k \text{ even}.
\end{cases}
\end{equation}

\eqref{Jla} is proved by
 standard contour integration arguments and implies  asymptotic
expansions for integrals $\int e^{i\lambda s^{k}} \chi(s) ds$ with $\chi\in
C^{\infty}_{0}$ (see \textit{e.g.} \S VIII.1.3 in \cite{stein}, or \S 7.7 in
\cite{Hoer}).

%\begin{lemma}\label{sklemma}
%Let $\eta\in C^2(\bbR)$ with compact support and assume that
%$\|\eta\|_{C^2}\le A$.
%Define for   $\lambda >1$
%\begin{equation}\label{standardasymptotics}
%J_\la= \int \eta(s) e^{i\la s^k}  ds
%\end{equation}
%Let \begin{equation}\label{ak}
%a_k=\begin{cases} \tfrac 2k\Gamma(\tfrac 1k)\sin(\tfrac{(k-1)\pi}{2k}),
%\, &k \text{ odd},\\
%\tfrac 2k\Gamma(\tfrac 1k)
%%e^{\pi i/(2k)}, \quad &k \text{ even}.
%\exp(i\tfrac{\pi}{2k}), \,  &k \text{ even}.
%\end{cases}
%\end{equation}
%Then  there is an absolute constant $C$ so that
%$$
%|J_\la-\eta(0) a_k \la^{-1/k}|
%\le C A\la^{-2/k}.
%$$
%\end{lemma}
We need small perturbations of such results. In what follows we set
$\|g\|_{C^{m}(I)}:=\max_{0\le j\le m}\sup_{x\in{I}} |g^{(j)}(x)|$.

\begin{lemma}
\label{asymptotics} Let $0<h\le 1$, $I=[-h,h]$, $I^*=[-2h,2h]$
 and let  $g\in C^{2}(I^*)$.
Suppose that
\begin{equation} \label{hassumption}
h\le \frac{1}{10(1+\|g\|_{C^2(I^*)})}
\end{equation}
and
%and suppose
%that $\sum_{j=0}^2\|g^{(j)}\|_\infty\le 1/10$.
let $\eta\in C^{1}$ be supported in $I$ and satisfy the bounds
\begin{equation}\label{etaassumption}
\|\eta\|_{\infty}+\|\eta^{\prime}\|_{1} \le A_{0},
\text{ and } \|\eta^{\prime}\|_{\infty}\le A_{1}.
\end{equation}
Let $k\ge 2$ and define
\begin{equation}
I_{\lambda}(\eta,x)= \int\eta(s) \exp\big(i\lambda(\sum_{j=1}^{k-2} x_{j}
s^{j}+ s^{k} + g(s) s^{k+1})\big) ds
\end{equation}
Let $\alpha_{k}$ be as in \eqref{ak}. Suppose $|x_{j}|\le\delta\lambda
^{(j-k)/k}$, $j=1,\dots, k-2$.
%, and $1/(2k!)<x_k\le2$.
Then there is an absolute constant $C$ so that, for $\lambda>2$,
\[
|I_{\lambda}(\eta,x)-\eta(0) \alpha_{k} \lambda^{-1/k}| \le C[A_{0}%
\delta\lambda^{-1/k}+A_{1}\lambda^{-2/k}(1+\beta_{k} \log\lambda)];
\]
here $\beta_{2}=1$, and $\beta_{k}=0$ for $k>2$.
\end{lemma}

\begin{proof}
%In view of our assumptions $\la x_k>1$ we can renormalize
%and  assume that $x_k=1$, after replacing $\lambda$ with $\lambda x_k$.

We set $u(s):=s(1+sg(s))^{1/k}$; then
\[u'(s)=(1+sg(s))^{-1+1/k}(1+sg(s)+k^{-1}s)\] and by our assumption on $g$ we quickly
verify that $(9/10)^{1/k}\le u^{\prime}(s)\le (11/10)^{1/k}$ for $-h\le s\le h$.
Thus  $u$ defines a
valid change of variable, with $u(0)=0$ and $u^{\prime}(0)=1$.
Denoting the
inverse by $s(u)$ we get
\[
I_{\lambda}(\eta,x)= \int\eta_{1}(u) \exp\big(i\lambda(\sum_{j=1}^{k-2} x_{j}
s(u)^{j}+ u^{k})\big)  du
\]
with $\eta_{1}(u)=\eta(s(u)) s^{\prime}(u)$.
Clearly $\eta_1$ is supported in $(-2h,2h).$
We observe that
\begin{equation} \label{etaone}
\|\eta_1\|_\infty+\|\eta_1'\|_1\lc A_0, \quad
\text{ and } \|\eta_1'\|_\infty \lc(A_0 h^{-1}+A_1).
\end{equation}
%Indeed, by differentiation of $u=s(u)(1+s(u)g(s(u)))^{1/k}$
%we get that
%%\begin{multline*}
%\[
%s'(u)=\frac{1}{L(s(u))}, \quad s''(u)=-\frac{L'(s(u))}{L(s(u))^3}
%\]
%where
%%\\\text{ where }
%$L(s)=(1+sg(s))^{(1-k)/k}(1+(1+k^{-1})sg(s)+k^{-1}s^2 g'(s)).$
Indeed implicit differentiation
and use of the assumption
\eqref{hassumption}
 reveals that  $|s''(u)|\lc (1+\|g\|_\infty)\lc h^{-1}$.
 Taking into account the support properties of
$\eta_1$ we obtain  \eqref{etaone}.

 In order to estimate certain
error terms we shall introduce dyadic decompositions. Let  $\chi_0\in C^\infty_0(\bbR)$ so that
\begin{equation} \label{chi0def}
\chi_{0}(s)=\begin{cases}1, \,\text{ if }|s|\le1/4,\\
0, \,\text{ if }|s|\ge1/2,
\end{cases}
\end{equation}
and
$m\ge1$,  define
\begin{equation}  \label{chimdef}
\chi_{m}(s)= \chi_{0}(2^{-m}s)-\chi_{0}(2^{-m+1}s),
\qquad m\ge 1.
\end{equation}

%We let $\zeta$ be a $C^\infty_0$ function satisfying $\zeta(u)=1$
%on a bounded interval containing the support of $\eta$.
We now split
\[
I_{\lambda}(\eta,x)=\eta_{1}(0) J_{\lambda}+ \sum_{m\geq0} E_{\lambda,m}+
\sum_{m\geq0} F_{\lambda,m}(x)
\]
where $J_{\lambda}$ is defined in \eqref{Jla} and
\begin{align*}
E_{\lambda,m}  &  = \int(\eta_{1}(u)-\eta_{1}(0)) \chi_{m}(\lambda^{1/k}u)
e^{i\lambda u^{k}} du,\\
F_{\lambda,m}(x)  &  = \int\eta_{1}(u) \Big( \exp(i\lambda(\sum_{j=1}^{k-2}
x_{j} s(u)^{j}))-1\Big) \chi_{m}(\lambda^{1/k}u) e^{i\lambda u^{k}} du.
\end{align*}
%with
%\begin{align*}\eta_2(u)&= \zeta(u) \int_0^1 \eta_1'(\sigma u)
%d\sigma\\
%\eta_3(u,x')&=\eta_1(u) \Big(
%\exp(i\la (\sum_{j=1}^{k-2} x_j s(u)^j))-1\Big).
%\end{align*}

In view of \eqref{Jla} the main term in our asymptotics is contributed by
$\eta_{1}(0) J_{\lambda}$ since $\eta_{1}(0)=\eta(0)$.

Now we estimate the terms $E_{\lambda,m}$. It is immediate that from an
estimate using the support of the amplitude that
\[
|E_{\lambda,0}|\le C\|\eta_{1}^{\prime}\|_{\infty}\lambda^{-2/k}.
\]
For $m\ge1$ we integrate by parts once to get
\[
E_{\lambda,m}= \frac{i}{k\lambda}\int\frac{d}{du}\big[
(\eta_{1}(u)-\eta_{1}(0))u^{1-k} \chi_{m}(\lambda^{1/k}u)\big]  e^{i\lambda
u^{k}} du
\]
and straightforward estimation gives
\[
|E_{\lambda,m}|\le C%
\begin{cases}
\|\eta_{1}^{\prime}\|_{\infty}2^{m(2-k)}\lambda^{-2/k}, \, & \text{ if }
2^{m}\le\lambda^{1/k}\\
\|\eta_{1}\|_{\infty}2^{m(1-k)}\lambda^{-1/k}, \, & \text{ if } 2^{m}>
\lambda^{1/k}.
\end{cases}
\]
Thus
\[
\sum_{m}|E_{\lambda,m}|\le C(\|\eta_{1}\|_{\infty}\lambda^{-2/k}(1+\beta_{k}
\log\lambda).
\]
%(and the logarithm is not present when $k>2$).

We now show that
\begin{equation}
\label{ila2}\sum_{m\geq0} |F_{\lambda,m}(x)|\le C[\|\eta_{1}\|_{\infty}%
+\|\eta^{\prime}\|_{1}] \delta\lambda^{-1/k}%
\end{equation}
and notice that only terms with $2^{m}\lambda^{-1/k}\le C$ occur in the sum.

Set $\zeta_{\lambda,x}(u)=\big(
\exp(i\lambda(\sum_{j=1}^{k-2} x_{j} s(u)^{j}))-1\big)$. For the term
$E_{\lambda,0}(x)$ we simply use the straightforward bound on the support of
$\chi_{0}(\lambda^{1/k}\cdot)$ which is (in view of $|s(u)|\approx|u|$)
\[
|\zeta_{\lambda,x}(u)|\le C\lambda\sum_{j=1}^{k-2} |x_{j}|| \lambda^{-j/k}|
\]
and since $|x_{j}|\le\delta\lambda^{-(k-j)/k}$ we get after integrating in
$u$
\[
|F_{\lambda,0}|\lesssim\|\eta_{1}\|_{\infty}\delta\lambda^{-1/k}.
\]

For $m>0$ we integrate by parts once and write
\begin{equation}
\label{intbypartsE}F_{\lambda,m}= i k^{-1}\lambda^{-1} \int\frac{d}{d
u}\big [u^{1-k} \chi_{m}(\lambda^{1/k} u) \eta_{1}(u) \zeta_{\lambda
,x}(u))\big] e^{i\lambda u^{k}} du.
\end{equation}
On the support of $\chi_{m}(\lambda^{1/k}\cdot)$,
\[
\begin{aligned} &| \zeta_{\la,x}(u)|\lesssim \lambda\sum_{j=1}^{k-2} \delta\lambda^{-(k-j)/k} (2^m\lambda^{-1/k})^j \lesssim \delta 2^{m(k-2)} \\ &\big|\tfrac{d}{du} \big[ u^{1-k} \chi_m(\lambda^{1/k}u)\big]\big|\lesssim \lambda 2^{-mk} \end{aligned}
\]
and also
\[
\begin{aligned} &|\zeta_{\la,x}'(u)| \lesssim \lambda\sum_{j=1}^{k-2} \delta\lambda^{-(k-j)/k} (2^m\lambda^{-1/k})^{j-1} \lesssim \delta 2^{m(k-3)}\lambda^{1/k} \\ &| u^{1-k} \chi_m(\lambda^{1/k}u)|\lesssim 2^{-m(k-1)}\lambda^{(k-1)/k}, \end{aligned}
\]
and thus we obtain the bound
\[
\int\Big|\frac{d}{du}\big [u^{1-k} \chi_{m}(\lambda^{1/k} u) \zeta_{\lambda
,x}(u)\eta_{1}(u))\big]\Big|du \lesssim[\|\eta_{1}\|_{\infty}+\|\eta
_{1}^{\prime}\|_{1}] \delta2^{-m}\lambda^{(k-1)/k}.
\]
Hence,
\[
\sum_{m\ge0}|F_{\lambda,m}| \lesssim\delta\lambda^{-1/k}
\]
which completes the proof of \eqref{ila2}.
\end{proof}

For the logarithmic lower bounds of $G_4(R)$ we shall need
some asymptotics for modifications of Airy functions.
Recall that for $t\in \bbR$ the Airy function is defined by the oscillatory integral
$$Ai(\tau) =\frac{1}{2\pi}\int_{-\infty}^{\infty} \exp(i(\tfrac {x^3}{ 3} +\tau x))
dx$$
and that for $t\to \infty$ we have
\begin{equation}
\label{airyasymp}
 Ai(-t)= \pi^{-1/2} t^{-1/4}
\cos(\tfrac 23 t^{3/2}-\tfrac \pi 4) (1+O(t^{-3/4})).
\end{equation}
This statement can be derived using the method of stationary phase
(combining expansions about the two critical points $\pm t^{1/2}$) or complex analysis arguments,
{\it cf.} \cite{erd} or  \cite{steinsh}, p. 330, see also an argument in
\cite{gsw}.

Let $g\in C^2([-1,1])$, and let $\eps>0$ be small, $\eps \ll
(1+\|g\|_{C^2})^{-1}$.
Let $\eta\in C^\infty_0$ with support in $(-\eps,\eps)$,  so that
$\eta(s)=1$ for $|s|\le \eps/2$.

\begin{lemma}\label{modairy}
%Suppose $\eta\in C^\infty_c$ satisfies $\eta(s)=1$ for $|s|<\eps$.
Define
\begin{equation}
J(\la, \vthb)=
\int e^{i\la (\tfrac {s^3} 3 -\vthb s)} e^{i\la g(s) s^4} \eta (s) ds
\end{equation}
Then, for  $0<\vthb<\eps^2/2$ and $\la>\eps^{-1}$
%\begin{equation}
\begin{align}
J(\la,\vthb) &= \la^{-1/3} Ai(-\la^{2/3} \vthb)+ E_1(\la, \vthb)
\notag
\\&= \pi^{-1/2} \la^{-1/2} \vthb^{-1/4} \cos \big(\tfrac 23 \la \vthb^{3/2}-
\tfrac \pi 4\big)+ E_2(\la, \vthb)
\label{Jlaasympt}
\end{align}
where, for $i=1,2$
\begin{equation}\label{errorJlaa}
|E_i(\la, \vthb)|\lc C_\eps\big[\la^{-1}\vthb^{-1}+\min\{\la \vthb^{5/2},
\vthb^{1/2}\}\big]
\end{equation}
\end{lemma}

\begin{proof}
We split
\begin{equation*}
J(\la,\vthb) = \sum_{i=1}^4J_i(\la,\vthb):=\sum_{i=1}^4
\int e^{i\la( \tfrac {s^3} 3 -\vthb s)} \zeta_i(s)\,ds
\end{equation*}
where
\begin{equation*}
%\begin{multline*}
\begin{aligned}
&\zeta_1(s)=1, \qquad \zeta_2(s)= (\eta(s)-1),\\
&\zeta_3(s) = \eta(s) (e^{i\la g(s) s^4}-1) \eta(C^{-1}\vthb^{-1/2}s),
\\
&\zeta_4(s) = \eta(s) (e^{i\la g(s) s^4}-1)(1-\eta(C^{-1}\vthb^{-1/2}s)).
\end{aligned}
\end{equation*}
%\end{multline*}
where $C\ge \eps^{-1}$.
By a scaling we see that
\[
J_1(\la,\vthb)=\la^{-1/3} Ai(-\la^{2/3} \vthb)
\] and we prove upper bounds for
the error terms $J_i$, $i=2,3,4$.
 Let \[\Phi(s)=-\vthb s+s^3/3\] then
$\Phi'(s)=-\vthb+s^2$ and
in the support of $\zeta_2$ we have $|\Phi'(s)|\ge c\eps$. Thus
 by an integration by  parts
$J_2(\la,\vthb )=O(\la^{-1})$.
Note that $\zeta_3$ is bounded and that also
$|\zeta_3(s)|\lc \la |\vthb|^2$. We integrate
over the support of $\zeta_3$ which is of length $O(\sqrt b)$ and  obtain
$J_3(\la,\vthb )=O(\min\{\vthb^{1/2},\la \vthb^{5/2}\})$.
To estimate $J_4(\la,\vthb)$ we
argue by van der Corput's Lemma, for the phases $\Phi$ and its
 perturbation  $\Psi(s):=\Phi(s)+s^4g(s)$.
Thus we split
\[
J_4(\la,\vthb) = \sum_m\sum_\pm J_{4, m,\pm}(\la,\vthb)
\]
where we have set
\[
J_{4,m,\pm}(\la,\vthb)=
\int e^{i\la\Psi(s)} \rho_{m,\pm} (s)\,ds
-
\int e^{i\la\Phi(s)} \rho_{m,\pm} (s)\,ds;
\]
here
$\rho_{m,+}(s)=\chi_{(0,\infty)}\eta(s) (1-\eta(C^{-1} \vthb^{-1/2}s))
\chi_m(C^{-1} \vthb^{-1/2}s)$, $\chi_m$ is as in
\eqref{chimdef} and $2^m \vthb^{1/2}\lc \eps$
(in view of the condition on $\eta$). Let
 $\rho_{m,-} $ is analogously defined, with support on $(-\infty,0)$.

We argue  as in the proof of Lemma
\ref{asymptotics}. Note that now $|\Phi'(s)|\approx 2^{2m}\vthb$,
$\partial_s(g(s) s^4)=O(2^{3m}\vthb^{3/2})$ and since
$2^m \vthb^{1/2}\lc \eps$ we also have
$|\Psi'(s)|\approx 2^{2m}\vthb$. Moreover observe
%that for $s\ge C\sqrt \vthb/2$
$\Phi''(s)= 2s+O(s^2)$
so that van der Corput's lemma can be applied can be
applied to the two integrals  defining $J_{4,m,\pm}(\la, \vthb)$.
We obtain
$J_{4,m,\pm}(\la, \vthb)=O(\la^{-1}\vthb^{-1} 2^{-2m})$.

Finally,
by \eqref{airyasymp} and \eqref{Jlaasympt},
the difference of $E_1$ and $E_2$ is
 $O(\la^{-1} \vthb^{-1})$.
This concludes the proof.
\end{proof}

\section{Lower bounds}

\label{frombelow}

For $w\in{\mathbb{R}}^{d}$ (usually restricted to the unit sphere), define
\begin{equation}
\label{cFRK}F_{R}(w)=\int\chi(t) e^{iR\langle\gamma(t),w\rangle} dt.
\end{equation}

The following result establishes inequality \eqref{withoutlog} of Theorem
\ref{lowerthm}.

\begin{proposition}
\label{sigmalowerthm} Suppose that for some $t_{0}\in I$ the vectors
$\gamma^{\prime}(t_{0})$, ..., $\gamma^{(k)}(t_{0}) $ are linearly
independent. Then $\chi\in C^{\infty}_{0}$ in \eqref{mudef} can be chosen so
that, for sufficiently large $R$,
%we have the lower bounds%
\begin{equation}
\label{sigmalowerbd}\|F_{R}\|_{L^{q}(S^{d-1})}\ge C R^{-\tfrac1k-\tfrac
{k^{2}-k-2}{2kq}}.
\end{equation}

\end{proposition}

\begin{proof}
We may assume $t_{0}=0$.
%We assume that for each $t$ the vectors
%$\gamma^{(j)}(t)$, $j=1,\dots, k$
%are linearly independent.
By a scaling and rotation we may assume that $\gamma^{(k)}(0)=e_{k}$. We shall
then show the lower bound $|F_{R}(\omega)|\ge c_{0} R^{-1/k}$ for a
neighborhood of $e_{k}$ which is of measure $\approx R^{-(k^{2}-k-2)/2k}$. Now
let $A_{k}$ be an invertible linear transformation which maps $e_{k}$ to
itself, and for $j=1,\dots, k-1$ maps $\gamma^{(j)}(0)$ to $e_{j}$,
$j=1,\dots, k$. Then the map $\omega\to(A_{k}^{*})^{-1}\omega/|(A_{k}
^{*})^{-1}\omega|$ defines a diffeomorphism from a spherical neighborhood of
$e_{k}$ to a spherical neighborhood of $e_{k}$. Thus we may assume for what
follows that $\gamma:[-1,1]\to{\mathbb{R}}^{d}$ satisfies
\begin{equation}
\gamma^{(j)}(0)=e_{j}, \qquad j=1,\dots, k.
\end{equation}
We may also assume that the cutoff function $\chi$ is supported in a small
open interval $(-\varepsilon, \varepsilon)$ so that $\chi(0)=1$.

As we have $\langle e_{k},\gamma^{(k-1)}(0)\rangle=0$ and $\langle
e_{k},\gamma^{k}(0)\rangle=1$ we can use the implicit function theorem to find
a neighborhood ${\mathcal{W}}_{k}$ of $e_{k}$ and an interval ${\mathcal{I}%
}_{k}=(-\varepsilon_{k},\varepsilon_{k})$ containing $0$ so that for all
$w\in{\mathcal{W}}_{k}$ the equation $\langle w,\gamma^{(k-1)}(t)\rangle=0$
has a unique solution $\widetilde t_{k}(w)\in{\mathcal{I}}_{k}$. This solution
is also homogeneous of degree $0$, i.e. $\widetilde t_{k}(sw)=\widetilde
t_{k}(w)$ for $s$ near $1$), and we have $\widetilde t_{k}(e_{k})=0$.

\begin{lemma}
\label{Ukdbd} There is $\varepsilon_{0}>0$, $R_{0}>1$, and $c>0$ so that for
all positive $\varepsilon<\varepsilon_{0}$ and all $R>R_{0}$ the following
holds.
%Let $U_{k,\delta}(R)$ be the set of all $\omega\in S^{d-1}$ for which
%$|\omega-e_k|\le \eps$ and
%$|\inn{\omega}{\gamma^{(j)}(\widetilde t_{k}(\omega))}|\le
%\delta R^{(j-k)/k}$, for $j=1,\dots, k-2$.
Let
\begin{multline*}
U_{k,\varepsilon}(R)=\big\{\omega\in S^{d-1}: |\omega-e_{k}|\le\varepsilon\\
\text{ and }|\langle\omega,\gamma^{(j)}(\widetilde t_{k}(\omega))\rangle
|\le\varepsilon R^{(j-k)/k}, \text{ for $j=1,\dots, k-2$.}\big\}
\end{multline*}
Then the spherical measure of $U_{k,\delta}(R)$ is at least $c\varepsilon
^{d-1} R^{-\frac{k^{2}-k-2}{2k}}$.
\end{lemma}

\begin{proof}
In a neighborhood of $e_{k}$ we parametrize the sphere by
\[
\omega(y)= (y_{1},\dots, y_{k-1},\sqrt{1-|y|^{2}}, y_{k},\dots, y_{d-1}).
\]
We introduce new coordinates $z_{1},\dots, z_{d-1}$ setting
\begin{equation}
\label{zjcoord}z_{j}={\mathfrak{z}}_{j}(y)=
\begin{cases}
\langle\omega(y),\gamma^{(j)}(\widetilde t_{k}(\omega(y)))\rangle, \quad &
j=1,\dots, k-2,\\
y_{j}, & j=k-1,\dots,d-1.
\end{cases}
\end{equation}
Then it is easy to see that ${\mathfrak{z}}$ defines a diffeomorphism between
small neighborhoods of the origin in ${\mathbb{R}}^{d-1}$; indeed the
derivative at the origin is the identity map.

The spherical measure of $U_{k,\varepsilon}(R)$ is comparable to the measure
of the set of $z\in\mathbb{R}^{d-1}$ satisfying $\left\vert z_{j}\right\vert
\leq\varepsilon R^{(j-k)/k}$, for $j=1,\dots,k-2$, and $\left\vert
z_{j}\right\vert \leq\varepsilon$ for $k-1\leq j\leq d-1$, and this set has
measure $\approx\varepsilon^{d-1}R^{-\frac{k^{2}-k-2}{2k}}$.
\end{proof}

We now verify that for sufficiently small $\varepsilon$ and
sufficiently large
$R$
\begin{equation}
|F_{R}(\omega)|\geq c_{0}R^{-1/k},\quad \omega\in U_{k,\varepsilon}(R),
\label{FRKlbd}%
\end{equation}
with some positive constant $c_{0}$; by Lemma \ref{Ukdbd}, this of course
implies the bound \eqref{sigmalowerbd}. To see \eqref{FRKlbd} we set
\begin{equation}
a_{j}(\omega)=\langle\omega,\gamma^{(j)}(\widetilde{t}_{k}(\omega))\rangle,
\label{ajom}%
\end{equation}
$s=t-\widetilde{t}_{k}(\omega)$ and expand%
\begin{equation}
\langle\omega,\gamma(t)\rangle-\langle\omega,\gamma(\widetilde{t}_{k}%
(\omega))\rangle\label{expand}\\
=\sum_{j=1}^{k-2}a_{j}(\omega)\frac{s^{j}}{j!}+a_{k}(\omega)\frac{s^{k}}%
{k!}+{\mathcal{E}}_{k}(\omega,s)s^{k+1},
\end{equation}
with ${\mathcal{E}}_{k}(\omega,s)=\int_{\sigma=0}^{1}\tfrac{(1-\sigma)^{k}%
}{k!}\langle\omega,\gamma^{(k+1)}(\widetilde{t}_{k}(\omega)+\sigma
s)\rangle\,d\sigma.$ If $\varepsilon$ is sufficiently small then we can apply
Lemma \ref{asymptotics} with $\omega\in U_{k,\varepsilon}(R)$, and the
choice $\lambda=R\langle\omega,\gamma^{(k)}(\widetilde{t}_{k}(\omega
))\rangle/k!$, and the lower bound \eqref{FRKlbd} follows.
\end{proof}

We now formulate bounds for $q\ge(k^{2}+k+2)/2$ for the case that
$\gamma^{(k+1)}\equiv0$ for some $k<d$; this of course implies that the curve
lies in a $k$-dimensional affine subspace.

\begin{proposition}
\label{sigmakplanethm}Suppose that $\gamma$ is a polynomial curve with
$\gamma^{(k+1)}\equiv0$ and suppose that for some $t_{0}\in I$ the vectors
$\gamma^{\prime}(t_{0})$, ..., $\gamma^{(k)}(t_{0})$ are linearly independent.
Then $\chi$ in \eqref{mudef} can be chosen so that for sufficiently large $R$
\begin{equation}
\Vert F_{R}\Vert_{L^{q}(S^{d-1})}\geq C%
\begin{cases}
R^{-k/q}[\log R]^{1/q}, & \quad q=\tfrac{k^{2}+k+2}{2},\\
R^{-k/q}, & \quad q>\tfrac{k^{2}+k+2}{2}.
\end{cases}
\label{secsigmalowerbd}%
\end{equation}

\end{proposition}

\begin{proof}
We first note that the assumption $\gamma^{(k+1)}\equiv0$ implies that the
curve is polynomial and for any fixed $t_{0}$ it stays in the affine subspace
through $\gamma(t_{0})$ which is generated by $\gamma^{(j)}(t_{0})$,
$j=1,\dots,k$. We shall prove a lower bound for $\widehat{\mu}$ in a
neighborhood of a vector $e\in S^{d-1}$ where $e$ is orthogonal to the vectors
$\gamma^{(j)}(t_{0})$. After a rotation we may assume that
\[
\gamma(t)=(\gamma_{1}(t),\dots,\gamma_{k}(t),0,\dots,0).
\]
For $\omega\in S^{d-1}$, we split accordingly $\omega=(\omega^{\prime},
\omega^{\prime\prime})$ with small  $\omega^{\prime}\in{\mathbb{R}}^{k}$,
namely
\[
|\omega^{\prime}|\approx2^{-l}
\]
where $1\ll R2^{-l}\ll R$. As before, we solve the first degree the equation
$\langle\gamma^{(k-1)}(t),\omega\rangle=0$ (observe that this is actually
independent of $\omega^{\prime\prime}$) with $t=\widetilde{t}_{k}%
(\omega^{\prime})$; now $\widetilde{t}_{k}$ is homogeneous of degree $0$ as a
function on ${\mathbb{R}}^{k}$. Then
\begin{align*}
&  e^{-i\langle\omega,\gamma(\widetilde{t}_{k}(\omega))\rangle}
F_{R}(\omega)\\
&  =\int\chi(\widetilde{t}_{k}(\omega^{\prime})+s)\exp
\Big(   \sum_{j=1}
^{k-2}\langle\omega,\gamma^{(j)}(\widetilde{t}_{k}(\omega^{\prime}
))\rangle\tfrac{s^{j}}{j!}+\langle\omega,\gamma^{(k)}(\widetilde{t}_{k}
(\omega^{\prime}))\rangle\tfrac{s^{k}}{k!} \Big)  ds
\end{align*}
If $k\geq3$, let $V_{k,l}(R)$ be the subset of the unit sphere in
${\mathbb{R}}^{k}$ which consists of those $\theta\in S^{k-1}$ which satisfy
the conditions
\[
|\langle\gamma^{(\nu)}(\widetilde{t}_{k}(\theta)),\theta\rangle|\leq
\varepsilon(R2^{-l})^{\tfrac{k-\nu}{k}},\quad\nu=1,\dots,k-2.
\]
Observe that the spherical measure of $V_{k,l}(R)$ (as a subset of $S^{k-1}$)
is $(R2^{-l})^{-(k^{2}-k-2)/(2k)}$, by Lemma \ref{Ukdbd}.
%(on the unit sphere $S^{k-1}$)
Now, if $k=2$, define
\[
U_{2,l}(R):=\{\omega=(\omega^{\prime},\omega^{\prime\prime})\in S^{d-1}%
:|\omega-e_{3}|\leq\delta,2^{-l}\leq|\omega^{\prime}|<2^{-l+1}\}.
\]
If $3\leq k<d$ let
%\begin{multline}%
\[
U_{k,l}(R):=\{\omega\in S^{d-1}:|\omega-e_{k+1}|\leq\delta,2^{-l}\leq
|\omega^{\prime}|<2^{-l+1},\tfrac{\omega^{\prime}}{|\omega^{\prime}|}\in
V_{k,l}(R)\}.
\]
We need a lower bound for the spherical measure (on $S^{d-1}$) of $U_{k,l}(R)$
and using polar coordinates in ${\mathbb{R}}^{k}$ we see that it is at least
\[
c\varepsilon^{d-1}2^{-lk}(R2^{-l})^{-\tfrac{k^{2}-k-2}{2k}}.
\]

If $\varepsilon$ is small we obtain a lower bound $c(R2^{-l})^{-1/k}$ on this
set; this follows from Lemma \ref{asymptotics} with $\lambda\approx
R2^{-l}$. Thus
\begin{align*}
\int_{U_{k,l}(R)}\big|F_{R}(\omega)|^{q}d\sigma(\omega)  &  \geq
c_{\varepsilon}(R2^{-l})^{-q/k}2^{-lk}(R2^{-l})^{-(k^{2}-k-2)/({2k})}\\
&  =c_{\varepsilon}R^{-q/k-(k^{2}-k-2)/(2k)}2^{l(q/k-(k^{2}+k+2)/(2k))}.
\end{align*}
As the sets $U_{k,l}(R)$ are disjoint in $l$ we may now sum in $l$ for
$CR^{-1}\leq2^{-l}\leq c$ for a large $C$ and a small $c$. Then we obtain that
$\sum_{l}\int_{U_{k,l}(R)}\big|F_{R}(\omega)|^{q}d\sigma(\omega)$ is bounded
below by $cR^{-q/k-(k^{2}-k-2)/(2kq)},$ if $q<(k^{2}+k+2)/2$; this yields the
bound that was already proved in Proposition \ref{sigmalowerthm}. If
$q>(k^{2}+k+2)/2$ then we get the lower bound $cR^{-k}$ and for the exponent
$q=(k^{2}+k+2)/2$ we obtain the lower bound $cR^{-k}\log R$. This yields \eqref{secsigmalowerbd}.
\end{proof}

\begin{proposition}
\label{qklowerbd} Suppose that $3\le k\le d$ and that for some $t_{0}\in I$
the vectors $\gamma^{\prime}(t_{0})$, ..., $\gamma^{(k)}(t_{0}) $ are linearly
independent. Then $\chi\in C^{\infty}_{0}$ in \eqref{mudef} can be chosen so
that for sufficiently large $R$
%we have the lower bounds%
\begin{equation}
\label{FRKqlbd}\|F_{R}\|_{L^{q}(S^{d-1})}\ge C R^{-(k-1)/q} [\log R]^{1/q},
\text{ if } q=q_{k-1}=\tfrac{k^{2}-k+2}{2}.
\end{equation}

\end{proposition}

\begin{proof}
We start with the same reductions as in the proof of Proposition
\ref{sigmalowerthm}, namely we may assume $t_{0}=0$ and $\gamma^{(j)}%
(0)=e_{j}$ for $j=1,\dots,k$; we shall then derive lower bounds for
$F_{R}(\omega)$ for $\omega$ near $e_{k}$. As before denote by $\widetilde
{t}_{k}(\omega)$ the solution $t$ of $\langle\gamma^{(k-1)}(t),\omega
\rangle=0$, for $\omega$ near $e_{k}$. We may use the expansion
\eqref{expand}. Define the polynomial approximation
\[
P_{k}(s,\omega)=P_{k}\left(  s\right)  =\sum_{j=1}^{k-2}a_{j}(\omega
)\frac{s^{j}}{j!}+a_{k}(\omega)\frac{s^{k}}{k!}.
\]
Note that $a_{k}(\omega)$ is near $1$ if $\omega$ is near $e_{k}$. In what
follows we shall only consider those $\omega$ with
\[
a_{k-2}(\omega)<0.
\]
In our analysis we need to distinguish between the cases $k=3$ and $k>3$.

\medskip

\noi\textit{The case $k=3$.}

We let for small $\delta$
\[\cU_{R,j}=\{\omega\in S^{d-1}: |\omega-e_3|\le \delta, \,
-2^{j+1}R^{-2/3}\le a_1(\omega)\le - 2^{j}R^{-2/3}
\}.
\]
We wish to use the asymptotics of
Lemma \ref{modairy},
with the parameters
 \[
\vthb=\vthb(\om)=\frac{-2a_1(\omega)}{a_3(\om)}\] and
%$\lambda= R\inn{\gamma'''(t_2(\omega))}{\omega}$
$\lambda= Ra_3(\om)/2 $ ($\approx R$)
to derive a lower bound on  a portion of $\cU_{R,j}$ whenever
$\la^{-2/3} \ll \vthb(\om) \ll \lambda^{-6/11}$; \textit{ i.e.}
\begin{equation} \label{j3range}
\delta^{-1}  \le 2^{j}\le \delta \lambda^{4/33}
\end{equation}
where $\delta$ is  small (but independent of large $\la$).

The range
%We observe  that in the range
\eqref{j3range}
is chosen so that the error terms in
\eqref{errorJlaa} (with $\la\approx R$) are
$\ll R^{-1/2}\vthb^{-1/4}$
if $\delta$  is small;
indeed the term $\lambda^{-1}\vartheta^{-1}$ is controlled by
$C\delta^{3/4}\la^{-1/2}\vartheta^{-1/4}$ in view of the
first inequality in \eqref{j3range}  and the term
%condition
%$\lambda^{-1}\vartheta^{-1}\lc
%\delta^{3/4} \la^{-1/2}\vartheta^{-1/4}$
%follows from the
%first inequality in \eqref{j3range}  and the condition
$\lambda\vartheta^{5/2}$ is bounded by
$ C\delta^{11/4} \la^{-1/2}\vartheta^{-1/4}$ because of
 the second restriction.
%The third error term in \eqref{errorJlaa}, $O(\vartheta)$
%is under control in the  larger range $2^j\ll \la^{2/15}$.
Since the main term
in \eqref{Jlaasympt} can be written as
\[(2/\pi)^{1/2} R^{-1/2} a_3(\om)^{-1/2} \vthb(\omega)^{-1/4}
\cos (\tfrac 13 Ra_3(\om) \vthb(\om)^{3/2}-\tfrac \pi 4)
\]
it
dominates the error terms in the range \eqref{j3range},   provided that we stay away from the zeroes of
 the cosine term.
%
%while for the main term we have  by \eqref{airyasymp}
%\begin{equation}\label{mainterm}
%\la^{-1/3} Ai(-\la^{2/3} a)= \pi^{-1/2}\la^{-1/2}a^{-1/4}
%\cos(\tfrac 23(\lambda a^{3/2}-\pi/4)) +O(\la a)^{-1}.
%\end{equation}
%Since $a>\la^{-2/3}\delta^{-1}$ the $O(\la a)^{-1}$ term
% is $\ll \la^{-1/2}a^{-1/4}$ if $\delta$ is small.
%
%We need to localize further to avoid the zeroes of the cosine term  in
%\eqref{mainterm}. Let
To achieve the necessary further localization we let, for positive integers $n$,
\[\cU_{R,j,n}=\{\omega\in \cU_{R,j}
\,:\,\big|\tfrac 13
Ra_3(\om) \vthb(\omega)^{3/2}-\tfrac {\pi}{ 4}-\pi n\big|<\tfrac \pi 4
\}.
\]
Let $j$ be in the range \eqref{j3range}. We use $|b^{3/2}-a^{3/2}|\approx
(\sqrt a+\sqrt b) |b-a|$ for $0<b,a \ll 1$. Since $\vthb(\omega)$ can be used as
one of the coordinates
on the unit sphere we see that
the spherical measure of
$\cU_{R,j,n}$ is $\gc \delta^2 R^{-2/3} 2^{-j/2}$ for the about
$ 2^{3j/2}$ values of $n$ for which $n\approx  2^{3j/2}$,
and on those disjoint sets $\cU_{R,j,n}$ the value of $F_R(\omega)$ is
$\ge c R^{-1/3} 2^{-j/4}$.

This implies that, for  $j$ as in \eqref{j3range},
\[
\text{meas}\big(\{\omega\in \cU_{R,j}:
|F_R(\omega)|\ge c_\delta R^{-1/3} 2^{-j/4}\}\big) \ge c_\delta' 2^{j}
R^{-2/3},
\]
and
thus
\[
\int_{\cU_{R,j}}|F_R(\omega)|^4 d\sigma(\omega)\gc  R^{-2}.
\]
Since the sets $\cU_{R,j}$ are disjoint
we may sum in $j$ over the range  \eqref{j3range}
and obtain the lower bound $\|F_R\|_4\gc R^{-1/2} (\log R)^{1/4}$
(with an implicit constant depending on $\delta$).

\medskip

\noi \textit{The case $k>3$.}
We try to follow in spirit the proof of the case for $k=3$.
Notice that
\[
P_{k}^{(k-2)}(s)=a_{k-2}(\omega)+a_{k}(\omega)s^{2}/2
\]
has then two real roots, one of them being
\[
s_{1}(\omega)=\Big(\frac{-2a_{k-2}(\omega)}{a_{k}(\omega)}\Big)^{1/2},
\]
the other one $s_{2}=-s_{1}$.
The idea is now to use, for suitable $\omega$,
an  asymptotic  expansion
for the part where $s$ is close to $s_1$, and, unlike in  the  case $k=3$,
 we shall  now be  able to neglect the contribution of the terms where $s$ is near $s_2$. To achieve this we define, for $j=1,\dots,k-3$,
\begin{equation}
\widetilde{a}_{j}(\omega)=P_{k}^{(j)}(s_{1}(\omega))=a_{j}(\omega)+\sum_{%
%TCIMACRO{\QATOP{1\leq\nu\leq k-2-j}{\QTR{group}{\text{or }\nu=k-j}}}%
%BeginExpansion
\genfrac{}{}{0pt}{}{1\leq\nu\leq k-2-j}{{\text{or }\nu=k-j}}%
%EndExpansion
}\frac{a_{j+\nu}(\omega)}{\nu!}\Big(\frac{-2a_{k-2}(\omega)}{a_{k}(\omega
)}\Big)^{\nu/2}.\label{tildeajom}
\end{equation}
%this definition is of course irrelevant when $k=3$.

We further restrict consideration to $\omega$ chosen in sets
%${\mathcal{V}}_{k,j}(\delta)$
%where for $k=3$
%\begin{equation}
%\label{V3j}{\mathcal{V}}_{3,j}(\delta)=\big\{\omega\in S^{d-1}: -2^{j+1}%
%R^{-2/k}< a_{1}(\omega)<-2^{j}R^{-2/k}, |e_{3}-\omega|\le\delta\big\}
%\end{equation}
%and for $k\ge4$
\begin{multline}
\label{Vkj}{\mathcal{V}}_{k,j}(\delta)=\big\{\omega\in S^{d-1}: -2^{j+1}%
R^{-2/k}< a_{k-2}(\omega)<-2^{j}R^{-2/k}, |e_{k}-\omega|\le\delta,\\
|\widetilde a_{\nu}(\omega)|\le\delta|a_{k-2}(\omega)|^{\frac{\nu}{2k-2}%
}R^{-\frac{k-\nu-1}{k-1}}, 1\le\nu\le k-3 \big\}.
\end{multline}

We shall see that if we choose $\omega$ from one of the sets ${\mathcal{V}%
}_{k,j}(\delta)$ with small $\delta$, and $j$ not too large then the main
contribution of the oscillatory integral comes from the part where
$|s-s_{1}(\omega)|\le s_{1}(\omega)/2$. We shall reduce to an application of
Lemma \ref{asymptotics} to derive a lower bound for that part. For the
remaining parts we shall derive smaller upper bounds using van der Corput's lemma.

For notational convenience we abbreviate
\[
b:=-a_{k-2}(\omega), \quad{\widetilde a}_{\nu}:=\widetilde a_{\nu}(\omega),
\quad s_{1}:=s_{1}(\omega), \quad{\widetilde t}_{k}:={\widetilde t}_{k}%
(\omega).
\]

%Fix a small $\eps>0$.
We now split
\begin{equation}
\label{IRER}e^{-i\langle\omega,\gamma({\widetilde t}_{k}(\omega))\rangle}%
F_{R}(\omega)= I_{R}(\omega) + E_{R}(\omega),
\end{equation}
where
\[
I_{R}(\omega) = \int\chi({\widetilde t}_{k}+s)\chi_{0}(20\frac{s-s_{1}}{s_{1}})
\exp(i R[ P_{k}(s)+s^{k+1}{\mathcal{E}}_{k+1} (s,\omega)])ds.
\]
Here $\chi_0$ is as in \eqref{chi0def} and thus the integrand is supported where
$|s-s_1|\le s_1/40$.

%The main term will be $I_{R,0}$ but we first derive upper bounds for the
%terms $I_{R,m}$, $m\ge 1$.

Notice that $P_{k}^{(k-1)}(s_{1})=s_{1}a_{k}(\omega)$ and $P_{k}%
^{(k)}(s)\equiv a_{k}(\omega)$. Let
\[
Q_{k-1}(s)=\sum_{\nu=1}^{k-3}\widetilde{a}_{\nu}\frac{(s-s_{1})^{\nu}}{\nu
!}+a_{k}s_{1}\frac{(s-s_{1})^{(k-1)}}{(k-1)!};
\]
then $P_{k}(s)-P_{k}(s_{1})=Q_{k-1}(s)+a_{k}(s-s_{1})^{k}/k!$.

Thus we can write
\[
I_{R}(\omega) = \int\eta(s) e^{iR (Q_{k-1}(s) +a_k(s-s_1)^k/k!)} ds
\]
with
\[\eta(s) =\chi({\widetilde t}_{k}+s) \chi_{0}( 10  s_{1}^{-1}(s-s_{1}))
\exp(i R s^{k+1}{\mathcal{E}}_{k+1} (s,\omega)).
\]
Note that by \eqref{chi0def} the function
$\eta$ is supported where $20  s_{1}^{-1}|s-s_{1}|\le 1/2$, i.e. in
$[s_1-h,s_1+h]$ with $h=s_1/40$.
 Clearly
$\|\eta\|_{\infty}=O(1)$, and since $s_{1}\approx\sqrt b$ it is
straightforward to check that
\begin{equation} \label{boundsforeta}
\|\eta^{\prime}\|_{\infty}+b^{-1/2} \|\eta^{\prime}\|_{1} \lc
(1+b^{-1/2}+
Rb^{k/2}),
\end{equation}
thus also
\begin{equation} \label{brestrict}
\|\eta\|_{\infty}+\|\eta^{\prime}\|_{1}\lc 1 \text{ if }
b\le R^{-2/{(k+1)}}.
\end{equation}
Moreover, if $g(s)=k^{-1}/ s_1$
%\frac{s-s_1}{k s_1}\]
then we can write
\begin{multline*}
RQ_{k-1}(s) +a_k\frac{(s-s_1)^k}{k!}
\\=\frac{R a_{k}s_{1}}{(k-1)!} \Big(\sum_{\nu=1}^{k-3} x_{\nu
}(s-s_{1})^{\nu}+(s-s_{1})^{k-1}+(s-s_1)^k g(s-s_1)\Big)
\end{multline*}
where
$|x_{\nu}|\lesssim b^{-1/2} |\widetilde a_{\nu}|$. The conditions
$|\widetilde a_{\nu}|\le\delta b^{\frac{\nu}{2k-2}}R^{-\frac{k-\nu-1}{k-1}}$
imply that
\[
|x_{\nu}|\lesssim\delta(Rb^{1/2})^{-\frac{k-\nu-1}{k-1}}.
\]
We of course  have
$\|g\|_{C^2([-h,h])}\le s_1^{-1}$ on $I^{*}=[-s_1/10, s_1 /10]$; thus
$h=s_1/40\le 10^{-1}(1+\|g\|_{C^2})^{-1}$.

Changing variables $\tilde s=s-s_1$  puts us in the position to apply
Lemma \ref{asymptotics} for
perturbations of the phase $\tilde s\mapsto\lambda \tilde s^{k-1}$, with
$\lambda:=R|a_{k}|
s_{1} = R\sqrt{2a_{k} b}\approx Rb^{1/2} $, and we have the bounds
$A_0\le C$ (if $b\le R^{-2/(k+1)}$
and $A_1\le (1+Rb^{k/2})$ for the parameters in
Lemma \ref{asymptotics}.
We thus obtain  (\textit{cf.} \eqref{ak})
\begin{multline}
\label{mainlowerbd}\big|I_{R}(\omega) -\alpha_{k-1}\chi(s_{1}(\omega))
(R\sqrt{2a_{k}b})^{-1/(k-1)}\big|\\
\lesssim\delta(Rb^{1/2})^{-1/(k-1)}
+ b^{-1/2} (Rb^{1/2})^{-2/(k-1)}
\log(Rb^{1/2}),
\end{multline}
provided that $b\le R^{-2/(k+1)}\ll 1$.
We wish to use this lower bound on the sets
${\mathcal{V}}_{k,j}(\delta)$.
In order to efficiently apply \eqref{mainlowerbd}
%we will have to restrict
%the values of $j$ (linked with $b$ as $b\approx R^{-2/k}2^j$) so that
%$(Rb^{1/2})^{-1/(k-1)} \gg (Rb^{k/2}+b^{-1/2})
%(Rb^{1/2})^{-2/(k-1)+\eps}$; {\it i.e.}
%when
%$b\gg R^{-2/k}$ and
%$b\ll R^{-\tfrac 2k(\tfrac{k^2-2k}{k^2-k-1})}$.
%The second condition is less restrictive than the condition
%$b\le R^{-2/(k+1)}$ already imposed in
%\eqref{mainlowerbd}, indeed
%$\frac{2}{k+1}-\frac{2(k^2-2k)}{k(k^2-k-1)}=\frac{2}{(k+1)(k^2-k-1)}>0.$
we shall  choose
$j$ so that
\begin{equation}
\label{jrange}
 R^{-\tau_{1}+2/k} \le2^{j}\le
R^{-\tau_{2}+2/k}
\end{equation}
%This will be  accomplished  when $j$ is chosen so that
%$R^{-\tau_{1}+2/k} \le2^{j}\le
%R^{-\tau_{2}+2/k} $
with $\tau_{1}, \tau_{2}$ satisfying
\[
\frac2{k}>\tau_{1}>\tau_{2}>
\frac{2}{k+1}.
%\frac{2(k^2-2k)}{k(k^2-k-1)},
\]
so that the main term in \eqref{mainlowerbd}
dominates the error terms.

We now need to bound from below  the measure of the set
${\mathcal{V}}_{k,j}(\delta)$.
We use the coordinates \eqref{zjcoord} on the sphere in a
neighborhood of $e_{k}$. In view of the linear independence of $\gamma
^{\prime\prime},\dots, \gamma^{(k-1)}$ we can use the functions $a_{j}%
({\mathfrak{z}}(y))$, $j\in\{1,\dots, k-2\}$, \textit{cf. \eqref{ajom}}, as a
set of partial coordinates.

We may also
%replace the coordinates
change coordinates
\[
(a_{1},\dots, a_{k-3}, a_{k-2}) \mapsto(\widetilde a_{1},\dots, \widetilde
a_{k-3}, a_{k-2}),
\]
with $a_{k-2}\equiv-b$; here we use the shear structure of the (nonsmooth)
change of variable \eqref{tildeajom}. Thus, as in Lemma \ref{Ukdbd}, we obtain
a lower bound for the spherical measure of ${\mathcal{V}}_{k,j}(\delta)$,
namely
%for $k=3$ \[
%\big|{\mathcal{V}}_{3,j}\big|\ge c \delta2^{j} R^{-2/3}%\]
%and, if $k>3$,
\begin{align*}
\big|{\mathcal{V}}_{k,j}(\delta)\big|  &  \ge c \delta^{d-2} 2^{j}R^{-\frac2k}
\prod_{\nu=1}^{k-3}\Big( \big(2^{j}R^{-\frac2k}\big)^{\frac{\nu}{2(k-1)}}
R^{-\frac{k-\nu-1}{k-1}}\Big)\\
&  =c \delta^{d-2} 2^{j} R^{-\frac2k} (2^{j} R^{-\frac2k})^{\frac
{(k-3)(k-2)}{4(k-1)}} R^{\frac{(k-3)(k-2)}{2(k-1)}-(k-3)}\\
&  =c \delta^{d-2} 2^{j\frac{k^{2}-k+2}{4(k-1)}} R^{-\frac{k^{2}-k-2}{2k}}%
\end{align*}
after a little arithmetic. Thus
%in both cases
\begin{equation}
\label{Vkjsize}\big|{\mathcal{V}}_{k,j}(\delta)\big|\\
\ge c \delta^{d-2}2^{j q_{k-1}/(2k-2)} R^{1/k - (k-1)/2}.
\end{equation}

Now if $\delta$ is chosen small and then fixed, and
$R$ is chosen large then \eqref{mainlowerbd} implies
the lower bound
\begin{equation}
\label{IRlowerbound}|I_{R}(\omega)| \ge c_{\delta}\big(R \sqrt{2^{j} R^{-2/k}%
})^{-1/(k-1)}= c_{\delta}2^{-j/(2k-2)} R^{-1/k}, \quad\omega\in{\mathcal{V}%
}_{k,j}(\delta),
\end{equation}
provided that $R^{-\tau_{1}+2/k}\le2^{j}\le
R^{-\tau_{2}+2/k}$.
We shall verify that for $j\ge0$
\begin{equation}
\label{ERupperbound}|E_{R}(\omega)|\lesssim%
%\begin{cases}
%%R^{-1/3}2^{-j} \quad & \text{if } k=3\\
R^{-1/k}\big(2^{-j/(k-2)}+2^{-3j/(2k-6)}\big), \quad
%& \text{if } k>3\end{cases}
\omega\in{\mathcal{V}}_{k,j}(\delta),
\end{equation}
and from \eqref{IRlowerbound} and \eqref{ERupperbound} it follows that
\[
|F_{R}(\omega)|\ge c_{\delta}2^{-j/(2k-2)} R^{-1/k},\quad\omega\in
{\mathcal{V}}_{k,j}(\delta)
\]
if  $R^{-\tau_{1}+2/k} \le2^{j}\le
R^{-\tau_{2}+2/k} $.
By \eqref {Vkjsize} this implies for the same
 range
a lower bound which is independent of $j$,
\begin{align*}
\int_{{\mathcal{V}}_{k,j}(\delta)} |F_{R}(\omega)|^{q_{k-1}} d\omega\ge
c_\delta
R^{-\frac{q_{k-1}}k- \frac{k^{2}-k-2}{2k}}= c_{\delta} R^{-(k-1)}.
\end{align*}
We sum in $j$, $R^{-\tau_{1}+2/k} \le2^{j}\le
R^{-\tau_{2}+2/k} $; this yields, for large $R$,
\[
\Big(\int_{\cup_{j}{\mathcal{V}}_{k,j}(\delta)} |F_{R}(\omega)|^{q_{k-1}}
d\omega\Big)^{1/q_{k-1}} \ge c_\delta^{\prime}R^{-(k-1)/q_{k-1} } \big(\log
R\big)^{1/q_{k-1}}
\]
which is the desired bound.

It remains to prove the upper bounds \eqref{ERupperbound}
for the error term $E_R$.
It is given by
\[
E_{R}(\omega) = e^{iRP_k(0)} \int\chi({\widetilde t}_{k}+s)(1-\chi_{0}(20\frac{s-s_{1}}{s_{1}}))e^{iR\phi(s)}
ds
\]
where
\[
\phi(s)    = P_{k}(s)-P_{k}(0)+s^{k+1}{\mathcal{E}}_{k+1}(s).
\]
We use a simple application of van der Corput's lemma. Write $\phi$ as
%\begin{align*}
\[
\phi(s)
%= P_{k}(s)-P_{k}(0)+s^{k+1}{\mathcal{E}}_{k+1}(s)\\
  =Q_{k-1}(s)+a_{k} (s-s_{1})^{k}/k! +s^{k+1}{\mathcal{E}}_{k+1}(s).
\]
%\end{align*}
and observe
\begin{align*}
\phi^{(k-2)}(s)  &  =a_{k} (s-s_{1})(s+s_{1})/2 +O(s^{3}),\\
\phi^{(k-3)}(s)  &  =\widetilde a_{k-3} +a_{k}\frac{(s-s_{1})^{2}}2
\big(\frac{2s_{1}}{3}+\frac s3\big)+O(s^{4}) .
\end{align*}
The integrand of the integral defining $E_R$ is supported
  where $|s-s_1|\ge s_1/80$, and $|s-s_1|\le c$ for small $c$.
We see that \[|\phi^{(k-2)}(s)|\ge c_{0} b\] if in addition
 $|s+s_{1}|\ge s_{1}/10.$

If  $|s+s_{1}|\le s_{1}/10$, this lower bound breaks down; however,
we have
then
\[
|\phi^{(k-3)}(s)|\ge c  b^{3/2} - |\widetilde a_{k-3}(\om)|.
\]
Now on $\cV_{k,j}(\delta)$ we have the restriction
\[|\widetilde a_{k-3}(\om)| \le
\delta b^{\frac{k-3}{2(k-1)}}R^{-\frac{2}{k-1}}\le \delta b^{3/2}
\]
where the last inequality is equivalent to the
imposed condition
$b\ge  R^{-2/k}$  (which holds when  $j\ge 2$).
Thus if $\delta $ is small
we have  $|\phi^{(k-3)}(s)|\approx b^{3/2} $ if $|s+s_1|\le s_1/10$.

We now split the integral into three parts
(using appropriate adapted cutoff functions), namely where
(i) $|s+s_1|\le s_1/10$, or (ii) $s+s_1\ge s_1/10$, or
(iii) $s+s_1\le -s_1/10$.
For parts (ii) and (iii) we can use
van der Corput's lemma
with
$k-2$ derivatives   and  see that the corresponding
 integrals are bounded by  $C(Rb)^{-1/(k-2)}$.
Similarly for part (i), if $k>4$  we can use van der Corput's lemma
with $(k-3)$ derivatives to see that the
 the corresponding
 integral is  bounded by
$C (Rb^{3/2})^{-1/(k-3)}$. The case $k=4$ requires a slightly
different argument
(as we do not
necessarily have adequate  monotonicity properties on $\phi'$),
however in the region (i) we now have  $\phi''(s)=O(b)$,
$|\phi'(s)|\gc b^{3/2}$
and integrating by parts once gives the required bound
$O( (Rb^{3/2})^{-1})$ also in this case.
Since $b\approx
R^{-2/k}2^{j}$, the upper bound \eqref{ERupperbound} follows.
\end{proof}

%\newpage

\end{document}